\documentclass[10pt,leqno]{article}

\usepackage{amsmath,amssymb,amsthm,mathrsfs,dsfont}

\usepackage[margin=3cm]{geometry} 

\usepackage{titlesec,hyperref}

\usepackage{color}

\usepackage{fancyhdr}
\titleformat{\subsection}{\it}{\thesubsection.\enspace}{1pt}{}

\newtheorem{theo}{Theorem}[section]
\newtheorem{lemm}[theo]{Lemma}
\newtheorem{defi}[theo]{Definition}
\newtheorem{coro}[theo]{Corollary}
\newtheorem{prop}[theo]{Proposition}
\newtheorem{rema}[theo]{Remark}
\numberwithin{equation}{section}

\allowdisplaybreaks 

\begin{document}
\title{Local well-posedness and blow-up criteria for a two-component Novikov system in the critical Besov space
\hspace{-4mm}
}

\author{Wei $\mbox{Luo}^1$\footnote{E-mail:  luowei23@mail2.sysu.edu.cn} \quad and\quad
 Zhaoyang $\mbox{Yin}^{1,2}$\footnote{E-mail: mcsyzy@mail.sysu.edu.cn}\\
 $^1\mbox{Department}$ of Mathematics,
Sun Yat-sen University,\\ Guangzhou, 510275, China\\
$^2\mbox{Faculty}$ of Information Technology,\\ Macau University of Science and Technology, Macau, China}
\date{}
\maketitle
\hrule

\begin{abstract}
In this paper we mainly investigate the Cauchy problem of a two-component Novikov system. We first prove the local well-posedness of the system in Besov spaces $B^{s-1}_{p,r}\times B^s_{p,r}$ with $p,r\in[1,\infty],~s>\max\{1+\frac{1}{p},\frac{3}{2}\}$ by using the Littlewood-Paley theory and transport equations theory. Then, by virtue of logarithmic interpolation inequalities and the Osgood lemma, we establish the local well-posedness of the system in the critical Besov space $B^{\frac{1}{2}}_{2,1}\times B^{\frac{3}{2}}_{2,1}$. Moreover, we present two blow-up criteria for the system by making use of the conservation laws.\\

\vspace*{5pt}
\noindent {\it 2010 Mathematics Subject Classification}: 35Q53 (35B30 35B44 35C07 35G25)

\vspace*{5pt}
\noindent{\it Keywords}: A two component Novikov system; Besov spaces; local well-posedness; blow-up criteria.
\end{abstract}

\vspace*{10pt}

\tableofcontents
\section{Introduction}
  In this paper we consider the Cauchy problem for the following two-component Novikov system:
  \begin{align}
\left\{
\begin{array}{ll}
\rho_{t}=\rho_{x}u^2+\rho uu_{x},\\[1ex]
m_{t}=3u_xum+u^2m_x-\rho(u\rho)_x ,\\[1ex]
m=u-u_{xx},\\[1ex]
\rho|_{t=0}=\rho_{0}, m|_{t=0}=m_{0}.\\[1ex]
\end{array}
\right.
\end{align}
This system was proposed by Popowicz in \cite{Popowicz}, in which the author verified that (1.1) has a Hamiltonian structure as follows \cite{Popowicz}:
 \begin{align*}
 \begin{pmatrix}
 \rho \\
 m
 \end{pmatrix}_t
 =\widehat{\mathcal{K}}
 \begin{pmatrix}
 \frac{\delta H}{\delta \rho}\\
 \frac{\delta H}{\delta m}
 \end{pmatrix},
 \end{align*}
 where $H=\frac{1}{2}\int_{\mathbb{R}}mu-\rho^2 dx$ and
 \begin{align*}
 \widehat{\mathcal{K}}= \begin{pmatrix}
 \rho^{-1}\partial\rho^2(\partial^3-4\partial_x)^{-1}\rho^2\partial\rho^{-1} & 3\rho^{-1}\partial\rho^2(\partial^3-4\partial_x)^{-1}m^{\frac{1}{3}}\partial m^{\frac{2}{3}}\\
 3m^{\frac{2}{3}}\partial(\partial^3-4\partial_x)^{-1}\rho^2\partial\rho^{-1} & -\rho\partial\rho+9m^{\frac{2}{3}}\partial m^{\frac{1}{3}}(\partial^3-4\partial_x)^{-1}m^{\frac{1}{3}}\partial m^{\frac{2}{3}}
 \end{pmatrix}.
 \end{align*}
Note that $G(x)\triangleq \frac{1}{2}e^{-|x|}$ is the kernel of $(1-\partial^2_x)^{-1}$. Then $G\ast m=u$ and the system (1.1) can be expressed as the following hyperbolic type:
  \begin{align}
\left\{
\begin{array}{ll}
\rho_{t}=\rho_{x}u^2+\rho uu_{x},\\[1ex]
u_{t}=u^2u_x+\partial_xG\ast(u^3+\frac{3}{2}uu^2_x-\frac{1}{2}u\rho^2)+\frac{1}{2}G\ast(u^3_x-u_x\rho^2) ,\\[1ex]
\rho|_{t=0}=\rho_{0}, u|_{t=0}=u_{0}.\\[1ex]
\end{array}
\right.
\end{align}
By setting $\rho=0$, the system (1.1) reduces to
\begin{align}
m_{t}=3u_xum+u^2m_x, ~~m=u-u_{xx},
\end{align}
which is nothing but the famous Novikov equation derived in \cite{Novikov}.  It was showed that (1.3) possesses a bi-Hamiltonian structure and an infinite sequence of conserved quantities in \cite{Hone}. Moreover, it admits exact peakon solutions $u(t,x)=\pm\sqrt{c}e^{|x-ct|}$ with $c>0$. It is worth mentioning that the peakons are solitons and present the characteristic singularity of greatest height and
largest amplitude, which arise as solutions to the free-boundary problem for incompressible Euler equations over a flat bed, cf. the discussions in \cite{Constantin2,Constantin.Escher4,Constantin.Escher5,Toland}. The above properties imply that peakons can be regarded as good
approximations to exact solutions of the governing
equations for water waves.\\
$~~~~~~$ The local well-posedness for (1.3) was studied in \cite{Wu.Yin2,Wu.Yin3,Wei.Yan,Wei.Yan2}. Concretely, for initial profile $u_0\in H^s(\mathbb{R})$ with $s>\frac{3}{2}$, it was shown in \cite{Wu.Yin2,Wei.Yan2} that (1.3) has a unique solution in $C([0,T);H^s(\mathbb{R}))$. Moveover, the local well-posedness for (1.3) in Besov spaces $C([0,T);B^s_{p,r}(\mathbb{R}))$ with $s>\max(\frac{3}{2},1+\frac{1}{p})$ was proved in \cite{Wu.Yin3,Wei.Yan}. The global existence of strong solutions were established in \cite{Wu.Yin2} under some sign conditions and the blow-up phenomena of the strong solutions were shown in \cite{Wei.Yan2}. The global weak solutions for (1.3) were studied in \cite{Laishaoyong,Wu.Yin}. \\
$~~~~~~$The early motivation to investigate the Novikov equation is that it can be regarded as a generalization of the well-known Camassa-Holm equation (CH) \cite{Camassa, Constantin.Lannes}:
\begin{align}
m_t+2m_xu+mu_x=0, ~~m=u-u_{xx}.
\end{align}
The most difference between the Novikov equation and the Camassa-Holm equation is that the former has cubic nonlinearity and the latter has quadratic nonlinearity.
\\
$~~~~~~$The Camassa-Holm equation was derived as a model for shallow water waves \cite{Camassa, Constantin.Lannes}. It has been investigated extensively because of its great physical significance in the past two decades. The CH equation has a bi-Hamiltonian structure \cite{Constantin-E,Fokas} and is completely integrable \cite{Camassa,Constantin-P,Constantin.McKean}. The solitary wave solutions of the CH equation were considered in \cite{Camassa,Camassa.Hyman}, where the authors showed that the CH equation possesses peakon solutions of the form $Ce^{-|x-Ct|}$. Constantin and Strauss verified that the peakon solutions of the CH equation are orbitally stable in \cite{Constantin.Strauss}. This means that the shape of soliton is stable. This equation attracted attention also in the context of the relevance of integrable equations to the modelling of tsunami waves, cf. the discussions in \cite{Constantin3,Constantin.Johnson,Lakshmanan,Segur1,Segur2}. Also, an aspect of considerable interest is the finite/infinite propagation speed associated
to the solutions of the equation, see e.g. \cite{Constantin.JMP,Henry}.\\
$~~~~~~$The local well-posedness for the CH equation was studied in \cite{Constantin.Escher,Constantin.Escher2,Danchin,Guillermo}. Concretely, for initial profiles $u_0\in H^s(\mathbb{R})$ with $s>\frac{3}{2}$, it was shown in \cite{Constantin.Escher,Constantin.Escher2,Guillermo} that the CH equation has a unique solution in $C([0,T);H^s(\mathbb{R}))$. Moveover, the local well-posedness for the CH equation in Besov spaces $C([0,T);B^s_{p,r}(\mathbb{R}))$ with $s>\max(\frac{3}{2},1+\frac{1}{p})$ was proved in \cite{Danchin}. The global existence of strong solutions were established in \cite{Constantin,Constantin.Escher,Constantin.Escher2} under some sign conditions and it was shown in \cite{Constantin,Constantin.Escher,Constantin.Escher2,Constantin.Escher3} that the solutions will blow up in finite time when the slope of initial data was bounded by a negative quantity. The global weak solutions for the CH equation were studied in \cite{Constantin.Molinet} and \cite{Xin.Z.P}. The global conservative and dissipative solutions of CH equation were presented in \cite{Bressan.Constantin} and  \cite{Bressan.Constantin2}, respectively.\\
 $~~~~~~$ The one popular two-component generalization of the Camassa-Holm equation is the following integrable two-component
Camassa-Holm shallow water system (2CH) \cite{Constantin.Ivanov}:
\begin{align}
\left\{
\begin{array}{ll}
m_t+um_x+2u_xm+\sigma\rho\rho_x=0, \\[1ex]
\rho_t+(u\rho)_x=0,
\end{array}
\right.
\end{align}
where $m=u-u_{xx}$ and $\sigma=\pm1$. Local well-posedness for (2CH) with the initial
data in Sobolev spaces and in Besov spaces were established in \cite{Constantin.Ivanov,Escher.Yin,GuiGuilong}. The blow-up phenomena and global existence of strong solutions to (2CH) in
Sobolev spaces were derived in \cite{Escher.Yin, Guan.Yin, GuiGuilong}. The existence
of global weak solutions for (2CH) with $\sigma=1$ was investigated in \cite{Guan.weak}.\\
$~~~~~~$The other one is the modified two-component Camassa-Holm
system (M2CH) \cite{Holm.Naraigh}:
\begin{align}
\left\{
\begin{array}{ll}
m_t+um_x+2u_xm+\sigma\rho\overline{\rho}_x=0, \\[1ex]
\rho_t+(u\rho)_x=0,
\end{array}
\right.
\end{align}
where $m=u-u_{xx}$, $\rho=(1-\partial^2_x)(\overline{\rho}-\overline{\rho}_0)$ and $\sigma=\pm1$. Local well-posedness for (M2CH) with the initial
data in Sobolev spaces and in Besov spaces were established in \cite{Guan.Karlsen} and \cite{Kai.Yin} respectively. Blow up phenomena of strong solution to (M2CH) were derived in \cite{Guan.Karlsen}. The existence
of global weak solutions for (M2CH) with $\sigma=1$ was investigated in \cite{Guan.weak.modified}.  The global conservative and dissipative solutions of (M2CH) equation was proposed in \cite{Tan.Yin} and  \cite{Tan.Yin2} respectively.\\
$~~~~~~$To our best knowledge, the Cauchy problem of (1.2) has not been studied yet. In this paper we first investigate the local well-posedness of (1.2) with initial data in Besov spaces $B^{s-1}_{p,r}\times B^s_{p,r}$ with $s>\max\{\frac{3}{2},1+\frac{1}{p}\}$. The main idea is based on the Littlewood-Paley theory and transport equations theory. Then, we study the local well-posedness in the critical space $B^{\frac{1}{2}}_{2,1}\times B^{\frac{3}{2}}_{2,1}$. We use the Friedrich's method to construct a sequence $(\rho_n,u_n)$ to approach the solution. It seems that one can't obtain  $(\rho_n,u_n)$ is a Cauchy sequence in $B^{-\frac{1}{2}}_{2,1}\times B^{\frac{1}{2}}_{2,1}$ as usual. By virtue of logarithmic interpolation inequalities and the Osgood lemma, we deduce that $(\rho_n,u_n)$ is a Cauchy sequence in $B^{-\frac{1}{2}}_{2,\infty}\times B^{\frac{1}{2}}_{2,\infty}$. Since the space $B^{-\frac{1}{2}}_{2,\infty}$ has low regularity, it follows that there are a lot of troubles when dealing with the nonlinear term. However, making good use of Bony's decomposition we can overcome these difficulties. Finally, we present two blow-up criteria with the help of the ordinary differential equation for the flow generated by $-u^2(t,x)$.\\
 $~~~~~~$ The paper is organized as follows. In Section 2 we introduce some preliminaries which will be used in sequel. In Section 3 we prove the local well-posedness of (1.2) by using Littlewood-Paley theory and transport equations theory. Moreover, by virtue of logarithmic interpolation inequalities and the Osgood lemma, we show the local well-posedness of (1.2) in the critical space. Section 4 is devoted to the study of two blow-up criteria for strong solutions to (1.2).
\section{Preliminaries}
In this section, we first recall the Littlewood-Paley decomposition and Besov spaces (for more details to see \cite{B.C.D}).
Let $\mathcal{C}$ be the annulus $\{\xi\in\mathbb{R}^{d}\big|\frac{3}{4}\leq|\xi|\leq\frac{8}{3}\}.$ There exist radial functions $\chi$ and $\varphi$, valued in the interval $[0,1]$, belonging respectively to $\mathcal{D}(B(0,\frac{4}{3}))$ and $\mathcal{D}(\mathcal{C})$, and such that $$\forall\xi\in\mathbb{R}^{d},~\chi(\xi)+\sum_{j\geq0}\varphi(2^{-j}\xi)=1,$$
$$|j-j'|\geq2\Rightarrow Supp ~\varphi(2^{-j}\xi)\cap Supp ~\varphi(2^{-j'}\xi)=\emptyset,$$
$$j\geq1\Rightarrow Supp ~\chi(\xi)\cap Supp ~\varphi(2^{-j'}\xi)=\emptyset.$$
  Define the set $\widetilde{\mathcal{C}}=B(0,\frac{2}{3})+\mathcal{C}$. Then we have
$$|j-j'|\geq5\Rightarrow 2^{j'}\widetilde{\mathcal{C}}\cap 2^{j}\mathcal{C}=\emptyset.$$
Further, we have $$\forall\xi\in\mathbb{R}^{d},~\frac{1}{2}\leq\chi^{2}(\xi)+\sum_{j\geq0}\varphi^{2}(2^{-j}\xi)\leq1,$$
   Denote $\mathcal{F}$ by the Fourier transform and $\mathcal{F}^{-1}$ by its inverse. From now on, we write $h=\mathcal{F}^{-1}\varphi$ and $\widetilde{h}=\mathcal{F}^{-1}\chi$.
The nonhomogeneous dyadic blocks $\Delta_{j}$ are defined by
$$\Delta_{j}u=0~~~ if~~~ j\leq-2,~~~\Delta_{-1}u=\chi(D)u=\int_{\mathbb{R}^{d}}\widetilde{h}(y)u(x-y)dy,$$
$$ and,~~~\Delta_{j}u=\varphi(2^{-j}D)u=2^{jd}\int_{\mathbb{R}^{d}}h(2^{j}y)u(x-y)dy ~~~if~~ j\geq0,$$
$$S_{j}u=\sum_{j'\leq j-1}\Delta_{j'}u.$$
The nonhomogeneous Besov spaces are denoted by $B^{s}_{p,r}(\mathbb{R}^d)$, i.e.,
$$B^{s}_{p,r}(\mathbb{R}^d)=\big\{u\in S'\big{|}\|u\|_{B^{s}_{p,r}(\mathbb{R}^d)}=(\sum_{j\geq-1}2^{rjs}\|\Delta_{j}u\|^{r}_{L^{p}(\mathbb{R}^d)})^{\frac{1}{r}}<\infty\big\}.$$
Next we introduce some useful lemmas and propositions about Besov spaces which will be used in the sequel.

\begin{prop}\label{2}
\cite{B.C.D} Let $1\leq p_{1} \leq p_{2} \leq \infty$ and $1\leq r_{1} \leq r_{2} \leq \infty$, and let $s$ be a real number. Then we have
$$B^{s}_{p_{1},r_{1}}(\mathbb{R}^d)\hookrightarrow B^{s-d(\frac{1}{p_{1}}-\frac{1}{p_{2}})}_{p_{2},r_{2}}(\mathbb{R}^d).$$
If $s>\frac{d}{p}~or ~s=\frac{d}{p},~r=1$, we then have $$B^{s}_{p,r}(\mathbb{R}^d)\hookrightarrow L^{\infty}(\mathbb{R}^d).$$
\end{prop}

\begin{defi}
\cite{B.C.D} The nonhomogeneous paraproduct of $v$ and $u$ is defined by
$$T_{u}v\triangleq\sum_{j}S_{j-1}u\Delta_{j}v.$$
The nonhomogeneous remainder of $v$ and $u$ is defined by
$$R(u,v)\triangleq\sum_{|k-j|\leq1}\Delta_{k}u\Delta_{j}v.$$
We have the following  Bony decomposition
$$uv=T_{u}v+R(u,v)+T_{v}u.$$
\end{defi}

\begin{prop}\label{p1}
\cite{B.C.D} For any couple of real numbers $(s, t)$ with $t$ negative and any $(p, r_{1}, r_{2})$
in $[1,\infty]^{3}$, there exists a constant $C$ such that:
$$\|T_{u}v\|_{B^{s}_{p,r}(\mathbb{R}^d)}\leq C\|u\|_{L^{\infty}(\mathbb{R}^d)}\|D^{k}v\|_{B^{s-k}_{p,r}(\mathbb{R}^d)},$$
$$\|T_{u}v\|_{B^{s+t}_{p,r}(\mathbb{R}^d)}\leq C\|u\|_{B^{t}_{\infty,r_{1}}(\mathbb{R}^d)}\|D^{k}v\|_{B^{s-k}_{p,r_{2}}(\mathbb{R}^d)},$$
where $r=\min\{1,\frac{1}{r_{1}}+\frac{1}{r_{2}}\}.$
\end{prop}

\begin{prop}\label{p2}
\cite{B.C.D} A constant $C$ exists which satisfies the following inequalities.
Let $(s_{1}, s_{2})$ be in $\mathbb{R}^2$ and $(p_{1}, p_{2}, r_{1}, r_{2})$ be in $[1,\infty]^{4}$. Assume that
$$\frac{1}{p}=\frac{1}{p_{1}}+\frac{1}{p_{2}}~~and~~\frac{1}{r}=\frac{1}{r_{1}}+\frac{1}{r_{2}} .$$
If $s_{1}+s_{2}>0$, then we have, for any $(u,v)$ in $B^{s_{1}}_{p_{1},r_{1}}\times B^{s_{2}}_{p_{2},r_{2}}$,
$$\|R(u,v)\|_{B^{s_{1}+s_{2}}_{p,r}(\mathbb{R}^d)}\leq \frac{C^{s_{1}+s_{2}+1}}{s_{1}+s_{2}}\|u\|_{B^{s_{1}}_{p_{1},r_{1}}(\mathbb{R}^d)}\|v\|_{B^{s_{2}}_{p_{2},r_{2}}(\mathbb{R}^d)}.$$
If $r=1$ and $s_{1}+s_{2}=0$, then we have, for any $(u,v)$ in $B^{s_{1}}_{p_{1},r_{1}}\times B^{s_{2}}_{p_{2},r_{2}}$,
$$\|R(u,v)\|_{B^{0}_{p,\infty}(\mathbb{R}^d)}\leq C\|u\|_{B^{s_{1}}_{p_{1},r_{1}}(\mathbb{R}^d)}\|v\|_{B^{s_{2}}_{p_{2},r_{2}}(\mathbb{R}^d)}.$$
\end{prop}

\begin{coro}\label{3}
\cite{B.C.D} For any positive real number $s$ and any $(p, r)$ in $[1,\infty]^{2}$, the
space $L^{\infty}(\mathbb{R}^d)\cap B^{s}_{p,r}(\mathbb{R}^d)$ is an algebra, and a constant $C$ exists such that
$$\|uv\|_{B^{s}_{p,r}(\mathbb{R}^d)}\leq C(\|u\|_{L^{\infty}(\mathbb{R}^d)}\|v\|_{B^{s}_{p,r}(\mathbb{R}^d)}+\|u\|_{B^{s}_{p,r}(\mathbb{R}^d)}\|v\|_{L^{\infty}(\mathbb{R}^d)}).$$
If $s>\frac{d}{p}$ or $s=\frac{d}{p},~r=1$, then we have
$$\|uv\|_{B^{s}_{p,r}(\mathbb{R}^d)}\leq C\|u\|_{B^{s}_{p,r}(\mathbb{R}^d)}\|v\|_{B^{s}_{p,r}(\mathbb{R}^d)}.$$
\end{coro}

\begin{lemm}\label{Rj.est}
\cite{B.C.D} Let $\sigma>0$, $1\leq r\leq \infty$ and $1\leq p\leq p_1\leq \infty$. Let $v$ be a vector field over $\mathbb{R}^d$. Define $R_j=[v\cdot \nabla, \Delta_j]f$. There exists a constant $C$ such that
$$\|(2^{j\sigma}\|R_j\|_{L^p(\mathbb{R}^d)})_j\|_{l^r}\leq C(\|\nabla v\|_{L^\infty(\mathbb{R}^d)}\|f\|_{B^\sigma_{p,r}(\mathbb{R}^d)}+\|\nabla f\|_{L^{p_2}(\mathbb{R}^d)}\|\nabla v\|_{B^{\sigma-1}_{p_1,r}(\mathbb{R}^d)}),$$
where $\frac{1}{p_2}=\frac{1}{p}-\frac{1}{p_1}$. Further, if $\sigma<1$ then
$$\|(2^{j\sigma}\|R_j\|_{L^p(\mathbb{R}^d)})_j\|_{l^r}\leq C\|\nabla v\|_{L^\infty(\mathbb{R}^d)}\|f\|_{B^\sigma_{p,r}(\mathbb{R}^d)}.$$
\end{lemm}

\begin{lemm}\label{Morse}
(Morse-type estimate, \cite{B.C.D,Danchin}) Let $\sigma>\max\{\frac{d}{p},\frac{d}{2}\}$ and $(p, r)$ in $[1,\infty]^{2}$. For any $a\in B^{\sigma-1}_{p,r}(\mathbb{R}^d)$ and $b\in B^{\sigma}_{p,r}(\mathbb{R}^d)$, there exists a constant $C$ such that
$$\|ab\|_{B^{\sigma-1}_{p,r}(\mathbb{R}^d)}\leq C\|a\|_{B^{\sigma-1}_{p,r}(\mathbb{R}^d)}\|b\|_{B^{\sigma}_{p,r}(\mathbb{R}^d)}.$$
\end{lemm}

The following two lemmas are crucial to study well-posedness in the critical space $B^{\frac{1}{2}}_{2,1}(\mathbb{R})$.
\begin{lemm}\label{Critical}
For any $a\in B^{-\frac{1}{2}}_{2,\infty}(\mathbb{R})$ and $b\in B^{\frac{1}{2}}_{2,1}(\mathbb{R})$, there exists a constant $C$ such that
$$\|ab\|_{B^{-\frac{1}{2}}_{2,\infty}(\mathbb{R})}\leq C\|a\|_{B^{-\frac{1}{2}}_{2,\infty}(\mathbb{R})}\|b\|_{B^{\frac{1}{2}}_{2,1}(\mathbb{R})}.$$
\begin{proof}
Using Bony's decomposition, we have $ab=T_ab+T_{b}a+R(a,b)$. By virtue of Proposition \ref{p1}, we deduce that
\begin{align*}
\|T_ab\|_{B^{-\frac{1}{2}}_{2,\infty}(\mathbb{R})}\leq C\|a\|_{B^{-1}_{\infty,\infty}(\mathbb{R})}\|b\|_{B^{\frac{1}{2}}_{2,1}(\mathbb{R})}\leq C\|a\|_{B^{-\frac{1}{2}}_{2,\infty}(\mathbb{R})}\|b\|_{B^{\frac{1}{2}}_{2,1}(\mathbb{R})}
\end{align*}
and
\begin{align*}
\|T_ba\|_{B^{-\frac{1}{2}}_{2,\infty}(\mathbb{R})}\leq C\|b\|_{L^\infty(\mathbb{R})}\|a\|_{B^{-\frac{1}{2}}_{2,\infty}(\mathbb{R})}\leq C\|a\|_{B^{-\frac{1}{2}}_{2,\infty}(\mathbb{R})}\|b\|_{B^{\frac{1}{2}}_{2,1}(\mathbb{R})}.
\end{align*}
Taking advantage of Proposition \ref{p2}, we infer that
\begin{align*}
\|R(a,b)\|_{B^{-\frac{1}{2}}_{2,\infty}(\mathbb{R})}\leq \|R(a,b)\|_{B^{0}_{1,\infty}(\mathbb{R})}\leq C\|a\|_{B^{-\frac{1}{2}}_{2,\infty}(\mathbb{R})}\|b\|_{B^{\frac{1}{2}}_{2,1}(\mathbb{R})}.
\end{align*}
\end{proof}
\end{lemm}

\begin{lemm}\label{com.est}(\cite{B.C.D})
Let $f$ be a smooth function on $\mathbb{R}^d$. Assume that $f$ is homogeneous of degree $m$ away from a neighborhood of $0$. Let $s\in\mathbb{R}$ and $(p,r)\in[1,\infty]$. There exists a constant $C$ such that
$$\|[T_a,f(D)]u\|_{B^{s-m+1}_{p,r}}\leq C\|\nabla a\|_{L^\infty}\|u\|_{B^s_{p,r}}.$$
\end{lemm}

\begin{lemm} \cite{Danchin2}\label{log.est}
For any $f\in B^{\frac{3}{2}}_{2,1}(\mathbb{R})$, there exists a constant $C$ such that
$$\|f\|_{B^{\frac{1}{2}}_{2,1}(\mathbb{R})}\leq C\|f\|_{B^{\frac{1}{2}}_{2,\infty}(\mathbb{R})}\ln\bigg(e+\frac{\|f\|_{B^{\frac{3}{2}}_{2,1}(\mathbb{R})}}{\|f\|_{B^{\frac{1}{2}}_{2,\infty}(\mathbb{R})}}\bigg).$$
\end{lemm}

\begin{rema}
\cite{B.C.D} Let $s\in\mathbb{R},1\leq p,r\leq\infty$, the following properties hold true. \\
$(i)~B^s_{p,r}(\mathbb{R}^d)$ is a Banach space and continuously embedding into $\mathcal{S}'(\mathbb{R}^d)$, where $\mathcal{S}'(\mathbb{R}^d)$ is the dual space of the Schwartz space $\mathcal{S}(\mathbb{R}^d)$. \\
$(ii)$ If $p,r<\infty$, then $\mathcal{S}(\mathbb{R}^d)$ is dense in $B^s_{p,r}(\mathbb{R}^d)$.\\
$(iii)$ If $u_n$ is a bounded sequence of $B^s_{p,r}(\mathbb{R}^d)$, then an element $u\in B^s_{p,r}(\mathbb{R}^d)$ and a subsequence $u_{n_k}$ exist such that
$$ \lim_{k\rightarrow\infty}u_{n_k}=u~~in~~\mathcal{S}'(\mathbb{R}^d)~~and~~\|u\|_{B^s_{p,r}(\mathbb{R}^d)}\leq C\liminf_{k\rightarrow\infty}\|u_{n_k}\|_{B^s_{p,r}(\mathbb{R}^d)}.$$
$(iv)~~B^s_{2,2}(\mathbb{R}^d)=H^s(\mathbb{R}^d)$.
\end{rema}

The following Osgood lemma appears as a substitution for Gronwall's lemma.
\begin{lemm}(Osgood's lemma, \cite{B.C.D})\label{Osgood}
Let $\rho\geq0$ be a measurable function, $\gamma>0$ be a locally integrable function and $\mu$ be a continuous and increasing function. Assume that, for some nonnegative real number $c$, the function $\rho$ satisfies
$$\rho(t)\leq c+\int^t_{t_0}\gamma(t')\mu(\rho(t'))dt'.$$
If $c>0$, then $-\mathcal{M}(\rho(t))+\mathcal{M}(c)\leq \displaystyle\int^t_{t_0}\gamma(t')dt'$ with $\mathcal{M}(x)=\displaystyle\int^1_x\frac{dr}{\mu(r)}$. \\ If $c=0$ and $\mu$ satisfies the condition $\displaystyle\int^1_0\frac{dr}{\mu(r)}=+\infty$, then the function $\rho=0$.
\end{lemm}

\begin{rema}
In this paper, we set $\mu(r)=r(1-\ln r)$ which satisfies the condition $\displaystyle\int^1_0\frac{dr}{\mu(r)}=+\infty$. A simple calculation shows that  $\mathcal{M}(x)=\ln(1-\ln x)$. Then, we deduce that
$$\rho(t)\leq \displaystyle c^{\bigg(\exp{\displaystyle\int^t_{t_0}-\gamma(t')dt'}\bigg)},~~~\textit{if}~~c>0.$$
\end{rema}

Now we introduce a priori estimates for the following transport equation
 \begin{align}
\left\{
\begin{array}{ll}
f_{t}+v\nabla f=g,\\[1ex]
f|_{t=0}=f_{0}.\\[1ex]
\end{array}
\right.
\end{align}
\begin{lemm}\label{est1}
(A priori estimates in Besov spaces, \cite{B.C.D}) Let $1\leq p \leq p_1\leq \infty$, $1\leq r\leq \infty$, $\sigma\geq -d\min(\frac{1}{p_1},\frac{1}{p'})$. For the solution $f\in L^{\infty}(0,T;B^\sigma_{p,r}(\mathbb{R}^d))$ of (2.1) with the velocity $\nabla v\in L^1(0,T;B^\sigma_{p,r}(\mathbb{R}^d)\cap L^{\infty}(\mathbb{R}^d))$,  the initial data $f_0\in B^\sigma_{p,r}(\mathbb{R}^d)$ and $g\in L^1(0,T;B^\sigma_{p,r}(\mathbb{R}^d))$, we have
\begin{align}
\|f(t)\|_{B^{\sigma}_{p,r}(\mathbb{R}^d)}\leq \|f_0\|_{B^\sigma_{p,r}(\mathbb{R}^d)}+\int^t_0\bigg(\|g(t')\|_{B^\sigma_{p,r}(\mathbb{R}^d)}+CV'_{p_1}(t')\|f(t')\|_{B^\sigma_{p,r}(\mathbb{R}^d)}\bigg)dt',\\
\|f\|_{L^{\infty}_t(B^{\sigma}_{p,r}(\mathbb{R}^d))}\leq \bigg(\|f_0\|_{B^\sigma_{p,r}(\mathbb{R}^d)}+\int^t_0\exp(-CV_{p_1}(t'))\|g(t')\|_{B^\sigma_{p,r}(\mathbb{R}^d)}dt'\bigg)\exp(CV_{p_1}(t)),
\end{align}
where  $V_{p_1}(t)=\displaystyle\int^t_0\|\nabla v\|_{B^{\frac{d}{p_1}}_{p_1,\infty}(\mathbb{R}^d)\cap L^{\infty}(\mathbb{R}^d)}$ if $\sigma<1+\frac{d}{p_1}$, $V_{p_1}(t)=\displaystyle\int^t_0\|\nabla v\|_{B^{\sigma-1}_{p_1,r}(\mathbb{R}^d)}$ if $\sigma>1+\frac{d}{p_1}$ or $\sigma=1+\frac{d}{p_1}, r=1$, and $C$ is a constant depending only on $\sigma,~p,~p_1$ and $r$.
\end{lemm}

\begin{lemm}\label{est2}
Let $1\leq p\leq \infty$, $1\leq r\leq \infty$, $\sigma> \max(\frac{1}{2},\frac{1}{p})$. For the solution $f\in L^{\infty}(0,T;B^\sigma_{p,r}(\mathbb{R}))$ of (2.1) with the velocity $ v\in L^1(0,T;B^{\sigma+1}_{p,r}(\mathbb{R}))$,  the initial data $f_0\in B^\sigma_{p,r}(\mathbb{R})$ and $g\in L^1(0,T;B^\sigma_{p,r}(\mathbb{R}^d))$, we have
\begin{align}
\|f\|_{L^{\infty}_t(B^{\sigma-1}_{p,r}(\mathbb{R}))}\leq \bigg(\|f_0\|_{B^{\sigma-1}_{p,r}(\mathbb{R})}+\int^t_0\exp(-CV(t'))\|g(t')\|_{B^{\sigma-1}_{p,r}(\mathbb{R})}dt'\bigg)\exp(CV(t)),
\end{align}
where  $V(t)=\displaystyle\int^t_0\|v\|_{B^{\sigma+1}_{p,r}(\mathbb{R})}$ and $C$ is a constant depending only on $\sigma,~p$ and $r$.
\begin{proof}
Applying $\Lambda=(1-\partial_{xx})^{-1}$ to both sides of (2.1) with $d=1$ yields that
\begin{align*}
(\Lambda f)_{t}+v(\Lambda f)_x=\Lambda g+ v(\Lambda f)_x-\partial_x\Lambda(vf)+\Lambda (v_x f).
\end{align*}
Taking advantage of Lemma \ref{est1}, we obtain that
\begin{align}
\|f\|_{B^{\sigma-1}_{p,r}}&\leq C\|\Lambda f\|_{B^{\sigma+1}_{p,r}}\leq C\exp(CV(t))\bigg\{\|\Lambda f_0\|_{B^{\sigma+1}_{p,r}}+\int^t_0\exp(-CV(t'))[\|\Lambda g(t')\|_{B^{\sigma+1}_{p,r}(\mathbb{R})}\\
\nonumber&+\|v(\Lambda f)_x-\partial_x\Lambda (v f)\|_{B^{\sigma+1}_{p,r}}+\|\Lambda(v_xf)\|_{B^{\sigma+1}_{p,r}}]dt'\bigg\}\\
\nonumber&\leq C\exp(CV(t))\bigg\{ \|f_0\|_{B^{\sigma-1}_{p,r}}+\int^t_0\exp(-CV(t'))[\| g(t')\|_{B^{\sigma-1}_{p,r}(\mathbb{R})}+\|v_xf\|_{B^{\sigma-1}_{p,r}}\\
\nonumber&+\|v(\Lambda f)_x-\partial_x\Lambda (v f)\|_{B^{\sigma+1}_{p,r}}]dt'\bigg\}.
\end{align}
By virtue of Lemma \ref{Morse}, we have 
\begin{align}
\|v_xf\|_{B^{\sigma-1}_{p,r}}\leq C \|v_x\|_{B^{\sigma}_{p,r}}\|f\|_{B^{\sigma-1}_{p,r}}\leq   C\|v\|_{B^{\sigma+1}_{p,r}}\|f\|_{B^{\sigma-1}_{p,r}}.
\end{align} 
Note that $v(\Lambda f)_x-\partial_x\Lambda (v f)=T_v(\Lambda f)_x+T_{(\Lambda f)_x}v+R(v,(\Lambda f)_x)-\partial_x\Lambda(T_vf)-\partial_x\Lambda(T_fv)-\partial_x\Lambda(R(f,v))$. In view of Proposition \ref{p1} and \ref{p2}, we have
\begin{align}
&\|T_{(\Lambda f)_x} v\|_{B^{\sigma+1}_{p,r}}\leq C\|(\Lambda f)_x\|_{L^\infty}\|v\|_{B^{\sigma+1}_{p,r}}\leq C\|f\|_{B^{\sigma-1}_{p,r}}\|v\|_{B^{\sigma+1}_{p,r}}, \\
&\|R(v,(\Lambda f)_x)\|_{B^{\sigma+1}_{p,r}}\leq C\|v\|_{B^1_{\infty,\infty}}\|(\Lambda f)_x\|_{B^\sigma_{p,r}}\leq C\|f\|_{B^{\sigma-1}_{p,r}}\|v\|_{B^{\sigma+1}_{p,r}}, \\
&\|\partial_x\Lambda(T_fv)\|_{B^{\sigma+1}_{p,r}}\leq C\|T_fv\|_{B^{\sigma-1}_{p,r}}\leq C\|f\|_{B^{-1}_{\infty,\infty}}\|v\|_{B^{\sigma+1}_{p,r}}\leq C\|f\|_{B^{\sigma-1}_{p,r}}\|v\|_{B^{\sigma+1}_{p,r}},\\
&\|\partial_x\Lambda(R(f,v))\|_{B^{\sigma+1}_{p,r}}\leq C\|R(f,v)\|_{B^\sigma_{p,r}}\leq C\|f\|_{B^{\sigma-1}_{p,r}}\|v\|_{B^1_{\infty,\infty}}\leq C\|f\|_{B^{\sigma-1}_{p,r}}\|v\|_{B^{\sigma+1}_{p,r}}.
\end{align} 
Since $T_v(\Lambda f)_x-\partial_x\Lambda(T_vf)=[T_v,\Lambda]f$, it follows from Lemma \ref{com.est} that
\begin{align}
\|T_v(\Lambda f)_x-\partial_x\Lambda(T_vf)\|_{B^{\sigma+1}_{p,r}}\leq C\|v_x\|_{L^\infty}\|f\|_{B^{\sigma-1}_{p,r}}\leq C\|f\|_{B^{\sigma-1}_{p,r}}\|v\|_{B^{\sigma+1}_{p,r}}.
\end{align}
Plugging (2.6)-(2.11) into (2.5), we deduce that
\begin{align}
\|f\|_{B^{\sigma-1}_{p,r}}\leq C\exp(CV(t))\bigg\{\|f_0\|_{B^{\sigma-1}_{p,r}}+\int^t_0\exp(-CV(t'))[\|g\|_{B^{\sigma-1}_{p,r}}+\|v\|_{B^{\sigma+1}_{p,r}}\|f\|_{B^{\sigma-1}_{p,r}}]dt'\bigg\},
\end{align}
or 
\begin{align}
\exp(-CV(t))\|f\|_{B^{\sigma-1}_{p,r}}\leq C\bigg\{\|f_0\|_{B^{\sigma-1}_{p,r}}+\int^t_0\exp(-CV(t'))[\|g\|_{B^{\sigma-1}_{p,r}}+\|v\|_{B^{\sigma+1}_{p,r}}\|f\|_{B^{\sigma-1}_{p,r}}]dt'\bigg\},
\end{align}
Applying Gronwall's inequality, we obtain the desire result.
\end{proof}
\end{lemm}

\begin{lemm}\label{est3}
 For the solution $f\in L^{\infty}(0,T;B^{1+\frac{1}{p}}_{p,r}(\mathbb{R}))$ of (2.1) with the velocity $ v\in L^1(0,T;B^{2+\frac{1}{p}}_{p,r}(\mathbb{R}))$,  the initial data $f_0\in B^{1+\frac{1}{p}}_{p,r}(\mathbb{R})$ and $g\in L^1(0,T;B^{1+\frac{1}{p}}_{p,r}(\mathbb{R}^d))$, we have
\begin{align}
\|f\|_{L^{\infty}_t(B^{1+\frac{1}{p}}_{p,r}(\mathbb{R}))}\leq \bigg(\|f_0\|_{B^{1+\frac{1}{p}}_{p,r}(\mathbb{R})}+\int^t_0\exp(-CV(t'))\|g(t')\|_{B^{1+\frac{1}{p}}_{p,r}(\mathbb{R})}dt'\bigg)\exp(CV(t)),
\end{align}
where  $V(t)=\displaystyle\int^t_0\|v\|_{B^{2+\frac{1}{p}}_{p,r}(\mathbb{R})}$ and $C$ is a constant depending only on $p$ and $r$.
\begin{proof}
Applying $\Lambda^{\frac{1}{2}}=(1-\partial_{xx})^{-\frac{1}{2}}$ to both sides of (2.1) with $d=1$ and by a similar method as in Lemma \ref{est2}, one can get the estimate. For the sake of conciseness, we omit the details here.
\end{proof}
\end{lemm}

{\bf Notations}. Since all function spaces in the following sections are over $\mathbb{R}$, for simplicity, we drop $\mathbb{R}$ in the notation of function spaces if there is no ambiguity.
\section{Local well-posedness}
In this section, we establish local well-posedness of the system (1.2) in Besov spaces. Our main results can be stated as follows.
\begin{theo}\label{th1}
Let $1\leq p,~r\leq \infty,~s>\max\{1+\frac{1}{p},~\frac{3}{2}\}$ and  $(\rho_0,~u_0)\in B^{s-1}_{p,r}\times B^s_{p,r}.$ There exists some $T>0$, such that the system (1.2) has a unique solution $(\rho,~u)$ in
   \begin{align}
   E^s_{p,r}(T)\triangleq
\left\{
\begin{array}{ll}
C([0,T);B^{s-1}_{p,r})\times C([0,T);B^s_{p,r})\cap C^1([0,T);B^{s-2}_{p,r})\times C^1([0,T);B^{s-1}_{p,r}) ,~~~~if~r<\infty, \\[1ex]
C_w([0,T);B^{s-1}_{p,\infty})\times C_w([0,T);B^s_{p,\infty})\cap C^{0,1}([0,T);B^{s-2}_{p,\infty})\times C^{0,1}([0,T);B^{s-1}_{p,\infty}),~~~~if~r=\infty.\\[1ex]
\end{array}
\right.
\end{align}
\end{theo}
Moreover, the solution map $S(t):(\rho_0,~u_0)\mapsto(S(t)\rho_0,~S(t)u_0)$ from $B^{s-1}_{p,r}\times B^s_{p,r}$ to $E^{s'}_{p,r}$ is H\"{o}lder continuous, that is
\begin{multline}
\|\rho_1-\rho_2\|_{L^\infty (0,T;B^{s'-1}_{p,r})}+\|u_1-u_2\|_{L^\infty (0,T;B^{s'}_{p,r})}\leq
C(\|\rho_1(0)-\rho_2(0)\|^{\theta}_{B^{s-1}_{p,r}}+\|u_1(0)-u_2(0)\|^{\theta}_{B^{s}_{p,r}}),
\end{multline}
where $\theta=s-s'\in(0,1]$.
\begin{rema}
Thanks to $B^s_{2,2}=H^s$, by taking $p=2, r=2$ in Theorem \ref{th1} one can get the local well-posedness of (1.2) in $C([0,T);H^{s-1})\times C([0,T);H^s)\cap C^1([0,T);H^{s-1})\times C^1([0,T);H^{s-1}),\quad s>\frac{3}{2}$.
\end{rema}

\begin{theo}\label{th2}
Let $(\rho_0,~u_0)\in B^{\frac{1}{2}}_{2,1}\times B^{\frac{3}{2}}_{2,1}.$ There exists some $T>0$, such that the system (1.2) has a unique solution $(\rho,~u)$ in
\begin{align*}
C([0,T);B^{\frac{1}{2}}_{2,1})\times C([0,T);B^{\frac{3}{2}}_{2,1})\cap C^1([0,T);B^{-\frac{1}{2}}_{2,1})\times C^1([0,T);B^{\frac{1}{2}}_{2,1}).
\end{align*}
Moreover,
\begin{multline}
\sup_{t\in[0,T)}\|\rho_{1}(t)-\rho_{2}(t)\|_{B^{s'-1}_{2,1}}+\sup_{t\in[0,T)}\|u_{1}(t)-u_{2}(t)\|_{B^{s'}_{2,1}}\\
\leq C (\|\rho_1(0)-\rho_2(0)\|_{B^{\frac{1}{2}}_{2,1}}+\|u_1(0)-u_2(0)\|_{B^{\frac{3}{2}}_{2,1}})^{\theta\exp\{-CT\}},
\end{multline}
where $\theta=\frac{3}{2}-s'\in(0,1]$.
\end{theo}
\begin{rema}
If $\rho=0$, then Theorem \ref{th1} and Theorem \ref{th2} cover the local well-posedness results of the Novikov equation obtained in \cite{Wu.Yin2,Wei.Yan2}.
\end{rema}
\subsection{Proof of Theorem \ref{th1}}
In order to prove Theorem \ref{th1}, we proceed as the following steps.\\

 {\bf Step 1:} First, we construct approximate solutions which are smooth solutions of some linear equations. Starting for $(\rho_1,u_1)\triangleq (S_{1}\rho_{0},S_{1}u_{0})$  we define by induction sequences $(\rho_{n},u_{n})_{n\geq1}$  by solving the following linear transport equations:
\begin{align}
\left\{
\begin{array}{ll}
\partial_t\rho_{n+1}-u^2_{n}\partial_x\rho_{n+1}=\rho_{n}u_n\partial_xu_n,\\[1ex]
\partial_tu_{n+1}-u^2_{n}\partial_xu_{n+1}=\partial_xG\ast(u^3_n+\frac{3}{2}u_n(\partial_xu_n)^2-\frac{1}{2}u_n\rho_n^2)+\frac{1}{2}G\ast((\partial_xu_n)^3-\partial_xu_n\rho^2_n),\\[1ex]
(\rho_{n+1},u_{n+1})|_{t=0}=(S_{n+1}\rho_{0},S_{n+1}u_{0}).\\[1ex]
\end{array}
\right.
\end{align}
We assume that $(\rho_n,u_n)\in L^{\infty}(0,T;B^{s-1}_{p,r})\times L^{\infty}(0,T;B^s_{p,r})$. Since $s>\max\{1+\frac{1}{p},2-\frac{1}{p}\}$, if follows that $B^{s-1}_{p,r}$ is an algebra, which leads to $\rho_{n}u_n\partial_xu_n\in L^{\infty}(0,T;B^{s-1}_{p,r})$ and $\partial_xG\ast(u^3_n+\frac{3}{2}u_n(\partial_xu_n)^2-u_n\rho_n^2)+\frac{1}{2}G\ast((\partial_xu_n)^3-\partial_xu_n\rho^2_n)\in L^{\infty}(0,T;B^s_{p,r})$. By the theory of transport equations, we obtain that $(\rho_{n+1}, u_{n+1})\in L^{\infty}(0,T;B^{s-1}_{p,r})\times L^{\infty}(0,T;B^s_{p,r})$.
For more details of the existence and uniqueness of the above system, one can refer to Chapter 3 in \cite{B.C.D}.\\

{\bf Step 2:} Next, we are going to find some positive $T$ such that for this fixed $T$ the approximate solutions are uniformly bounded on $[0,T]$.
We define that
$U_{n}(t)=\displaystyle\int^{t}_{0}\|u_n(t')\|^2_{B^{s}_{p,r}}dt'$.

If $s\neq 2+\frac{1}{p}$, by virtue of Lemma \ref{est1} with $p_1=p$, $\sigma=s$ and $p_1=p$, $\sigma=s-1$ respectively, we infer that
\begin{multline}
\|u_{n+1}\|_{L^{\infty}_t(B^{s}_{p,r})}\leq e^{CU_{n}(t)}\bigg(\|S_{n+1}u_0\|_{B^s_{p,r}}\\
+\int^t_0e^{-CU_{n}(t')}\|\partial_xG\ast(u^3_n+\frac{3}{2}u_n(\partial_xu_n)^2-\frac{1}{2}u_n\rho_n^2)+\frac{1}{2}G\ast((\partial_xu_n)^3-\partial_xu_n\rho^2_n)\|_{B^s_{p,r}}dt'\bigg),
\end{multline}
\begin{align}
\|\rho_{n+1}\|_{L^{\infty}_t(B^{s-1}_{p,r})}&\leq e^{CU_{n}(t)}\bigg(\|S_{n+1}\rho_0\|_{B^{s-1}_{p,r}}
+\int^t_0e^{-CU_{n}(t')}\|\rho_n u_n\partial_xu_n\|_{B^{s-1}_{p,r}}dt'\bigg)\\
\nonumber&\leq e^{CU_{n}(t)}\bigg(\|S_{n+1}\rho_0\|_{B^{s-1}_{p,r}}
+C\int^t_0e^{-CU_{n}(t')}\|\rho_n\|_{B^{s-1}_{p,r}}\|u_n\|^2_{B^s_{p,r}}dt'\bigg),
\end{align}
Note that the operator $(1-\partial_{xx})^{-1}$ is a multiplier of degree $-2$ and $G\ast f=(1-\partial_{xx})^{-1}f$, we deduce that
\begin{align}
\|\partial_xG\ast(u^3_n+\frac{3}{2}u_n(\partial_xu_n)^2-\frac{1}{2}u_n\rho_n^2)\|_{B^s_{p,r}}&\leq C(\|u^3_n\|_{B^{s-1}_{p,r}}+\|u_n(\partial_xu_n)^2\|_{B^{s-1}_{p,r}}+\|u_n\rho_n^2\|_{B^{s-1}_{p,r}})\\
\nonumber&\leq C(\|u_n\|^3_{B^{s-1}_{p,r}}+\|u_n\|_{B^{s-1}_{p,r}}\|\partial_xu_n\|^2_{B^{s-1}_{p,r}}+\|u_n\|_{B^{s-1}_{p,r}}\|\rho_n\|^2_{B^{s-1}_{p,r}})\\
\nonumber&\leq  C(\|u_n\|^3_{B^{s}_{p,r}}+\|u_n\|_{B^{s}_{p,r}}\|\rho_n\|^2_{B^{s-1}_{p,r}})
\end{align}
and
\begin{multline}
\|G\ast((\partial_xu_n)^3-\partial_xu_n\rho^2_n)\|_{B^s_{p,r}}\leq C(\|(\partial_xu_n)^3\|_{B^{s-2}_{p,r}}+\|\partial_xu_n\rho_n^2\|_{B^{s-2}_{p,r}})\\
\leq C(\|(\partial_xu_n)^3\|_{B^{s-1}_{p,r}}+\|\partial_xu_n\rho_n^2\|_{B^{s-1}_{p,r}})
\leq C(\|u_n\|^3_{B^{s}_{p,r}}+\|u_n\|_{B^{s}_{p,r}}\|\rho_n\|^2_{B^{s-1}_{p,r}}).
\end{multline}
Plugging (3.7) and (3.8) into (3.5), we obtain
\begin{multline}
\|u_{n+1}\|_{L^{\infty}_t(B^{s}_{p,r})}\leq e^{CU_{n}(t)}\bigg(\|S_{n+1}u_0\|_{B^s_{p,r}}
+C\int^t_0e^{-CU_{n}(t')}(\|u_n\|^3_{B^{s}_{p,r}}+\|u_n\|_{B^{s}_{p,r}}\|\rho_n\|^2_{B^{s-1}_{p,r}})dt'\bigg),
\end{multline}
which together with (3.6) and Young's inequality imply that
 \begin{multline}
\|\rho_{n+1}\|_{L^{\infty}_t(B^{s-1}_{p,r})}+\|u_{n+1}\|_{L^{\infty}_t(B^{s}_{p,r})}\leq e^{CU_{n}(t)}\bigg[C(\|\rho_0\|_{B^{s-1}_{p,r}}+\|u_0\|_{B^s_{p,r}})\\
+C\int^t_0e^{-CU_{n}(t')}(\|u_n\|_{B^{s}_{p,r}}+\|\rho_n\|_{B^{s-1}_{p,r}})^3dt'\bigg],
\end{multline}
where we take $C\geq 1$.
We fix a $T>0$ such that $4C^3T(\|\rho_0\|_{B^{s-1}_{p,r}}+\|u_0\|_{B^s_{p,r}})^2<1$ and suppose that
\begin{align}
\forall t\in[0,T],~~\|\rho_{n}\|_{B^{s-1}_{p,r}}+\|u_{n}\|_{B^{s}_{p,r}}\leq \frac{C(\|\rho_0\|_{B^{s-1}_{p,r}}+\|u_0\|_{B^s_{p,r}})}{\sqrt{1-4C^3(\|\rho_0\|_{B^{s-1}_{p,r}}+\|u_0\|_{B^s_{p,r}})^2 t}}\leq \frac{C(\|\rho_0\|_{B^{s-1}_{p,r}}+\|u_0\|_{B^s_{p,r}})}{\sqrt{1-4C^3(\|\rho_0\|_{B^{s-1}_{p,r}}+\|u_0\|_{B^s_{p,r}})^2 T}}.
\end{align}
Since $U_{n}(t)=\displaystyle\int^{t}_{0}\|u_n(t')\|^2_{B^{s}_{p,r}}dt'$, it follows that
\begin{align}
e^{CU_{n}(t)-CU_{n}(t')}\leq \exp{\bigg\{\int^t_{t'} \frac{C^3(\|\rho_0\|_{B^{s-1}_{p,r}}+\|u_0\|_{B^s_{p,r}})^2}{1-4C^3(\|\rho_0\|_{B^{s-1}_{p,r}}+\|u_0\|_{B^s_{p,r}})^2 t}d\tau\bigg\}}=\sqrt[4]{\frac{1-4C^3(\|\rho_0\|_{B^{s-1}_{p,r}}+\|u_0\|_{B^s_{p,r}})^2t'}{1-4C^3(\|\rho_0\|_{B^{s-1}_{p,r}}+\|u_0\|_{B^s_{p,r}})^2t}},
\end{align}
\begin{align}
e^{CU_{n}(t)}\leq \exp{\bigg\{\int^t_{0} \frac{C^3(\|\rho_0\|_{B^{s-1}_{p,r}}+\|u_0\|_{B^s_{p,r}})^2}{1-4C^3(\|\rho_0\|_{B^{s-1}_{p,r}}+\|u_0\|_{B^s_{p,r}})^2 t}d\tau\bigg\}}=\sqrt[4]{\frac{1}{1-4C^3(\|\rho_0\|_{B^{s-1}_{p,r}}+\|u_0\|_{B^s_{p,r}})^2t}}.
\end{align}
Plugging (3.11)-(3.13) into (3.10) yields that
 \begin{align}
&\|\rho_{n+1}\|_{L^{\infty}_t(B^{s-1}_{p,r})}+\|u_{n+1}\|_{L^{\infty}_t(B^{s}_{p,r})}\leq \\ \nonumber&\sqrt[4]{\frac{1}{1-4C^3(\|\rho_0\|_{B^{s-1}_{p,r}}+\|u_0\|_{B^s_{p,r}})^2t}}\bigg(\|\rho_0\|_{B^{s-1}_{p,r}}+\|u_0\|_{B^s_{p,r}}
+\displaystyle\int^t_0\frac{C^3(\|\rho_0\|_{B^{s-1}_{p,r}}+\|u_0\|_{B^s_{p,r}})^3}{\big[1-4C^3(\|\rho_0\|_{B^{s-1}_{p,r}}+\|u_0\|_{B^s_{p,r}})^2t'\big]^{\frac{5}{4}}}dt'\bigg)\\
\nonumber&\leq\frac{C(\|\rho_0\|_{B^{s-1}_{p,r}}+\|u_0\|_{B^s_{p,r}})}{\sqrt{1-4C^3(\|\rho_0\|_{B^{s-1}_{p,r}}+\|u_0\|_{B^s_{p,r}})^2 t}}.
\end{align}
Therefore, by induction, we obtain that $(\rho_n,u_n)$ is bounded in $L^{\infty}(0,T;B^{s-1}_{p,r})\times L^{\infty}(0,T;B^{s}_{p,r})$.

If $s=2+\frac{1}{p}$. Applying Lemma \ref{est3}, we obtain that
\begin{align*}
\|\rho_{n+1}\|_{L^{\infty}_t(B^{1+\frac{1}{p}}_{p,r})}&\leq e^{CU_{n}(t)}\bigg(\|S_{n+1}\rho_0\|_{B^{1+\frac{1}{p}}_{p,r}}
+\int^t_0e^{-CU_{n}(t')}\|\rho_n u_n\partial_xu_n\|_{B^{1+\frac{1}{p}}_{p,r}}dt'\bigg)\\
\nonumber&\leq e^{CU_{n}(t)}\bigg(\|S_{n+1}\rho_0\|_{B^{1+\frac{1}{p}}_{p,r}}
+C\int^t_0e^{-CU_{n}(t')}\|\rho_n\|_{B^{1+\frac{1}{p}}_{p,r}}\|u_n\|^2_{B^{2+\frac{1}{p}}_{p,r}}dt'\bigg),
\end{align*}
Taking advantage of Lemma \ref{est1} with $p_1=p$, $\sigma=2+\frac{1}{p}$ and by the similar argument as $s\neq 2+\frac{1}{p}$ , we deduce that
\begin{multline*}
\|u_{n+1}\|_{L^{\infty}_t(B^{2+\frac{1}{p}}_{p,r})}\leq e^{CU_{n}(t)}\bigg(\|S_{n+1}u_0\|_{B^{2+\frac{1}{p}}_{p,r}}
+C\int^t_0e^{-CU_{n}(t')}(\|u_n\|^3_{B^{2+\frac{1}{p}}_{p,r}}+\|u_n\|_{B^{2+\frac{1}{p}}_{p,r}}\|\rho_n\|^2_{B^{1+\frac{1}{p}}_{p,r}})dt'\bigg),
\end{multline*}
By the same token, we get
 $(\rho_n,u_n)$ is bounded in $L^{\infty}(0,T;B^{1+\frac{1}{p}}_{p,r})\times L^{\infty}(0,T;B^{2+\frac{1}{p}}_{p,r})$. \\

{\bf Step 3:} From now on, we are going to show that $\{(\rho_n,u_n)\}$ is a Cauchy sequence
 in some Banach space. For this purpose, we deduce from (3.4) that
  \begin{align}
\left\{
\begin{array}{ll}
\partial_t(\rho_{n+m+1}-\rho_{n+1})-u^2_{n+m}\partial_x(\rho_{n+m+1}-\rho_{n+1})=(u^2_{n+m}-u^2_{n})\partial_x\rho_{n+1}+R^1_{n,m},\\[1ex]
\partial_t(u_{n+m+1}-u_{n+1})-u^2_{n+m}\partial_x(u_{n+m+1}-u_{n+1})=(u^2_{n+m}-u^2_{n})\partial_xu_{n+1}+R^2_{n,m},\\[1ex]
\end{array}
\right.
\end{align}
where $R^1_{n,m}=\rho_{m+n}u_{m+n}\partial_xu_{n+m}-\rho_nu_n\partial_xu_n,~~ R^2_{n,m}=\partial_xG\ast[u^3_{m+n}-u^3_n+\frac{3}{2}[u_{m+n}(\partial_xu_{m+n})^2-u_n(\partial_xu_n)^2]+\frac{1}{2}(u_n\rho^2_n-u_{m+n}\rho^2_{m+n})]+\frac{1}{2}G\ast[(\partial_xu_{m+n})^3-(\partial_xu_n)^3+\partial_xu_n\rho^2_n-\partial_xu_{m+n}\rho^2_{m+n}].$\\
By virtue of Lemma \ref{est2} with $p_1=p$ and using the fact that $\{(\rho_n,u_n)\}$ is bounded in $L^{\infty}(0,T;B^{s-1}_{p,r})\times L^{\infty}(0,T;B^{s}_{p,r})$, we infer that
\begin{multline}
\|\rho_{n+m+1}-\rho_{n+1}\|_{L^{\infty}(0,T;B^{s-2}_{p,r})}\leq C\bigg(\|S_{n+m+1}\rho_0-S_{n+1}\rho_0\|_{B^{s-2}_{p,r}}\\
+\int^T_0\|(u^2_{n+m}-u^2_{n})\partial_x\rho_{n+1}\|_{B^{s-2}_{p,r}}+\|R^1_{n,m}\|_{B^{s-2}_{p,r}}dt'\bigg),
\end{multline}
\begin{multline}
\|u_{n+m+1}-u_{n+1}\|_{L^{\infty}(0,T;B^{s-1}_{p,r})}\leq C\bigg(\|S_{n+m+1}u_0-S_{n+1}u_0\|_{B^{s-1}_{p,r}}\\
+\int^T_0\|(u^2_{n+m}-u^2_{n})\partial_xu_{n+1}\|_{B^{s-1}_{p,r}}+\|R^2_{n,m}\|_{B^{s-1}_{p,r}}dt'\bigg).
\end{multline}
Taking advantage of Lemma \ref{Morse} with $\sigma=s-1$ and $d=1$, we have
\begin{multline}
\|(u^2_{n+m}-u^2_{n})\partial_x\rho_{n+1}\|_{B^{s-2}_{p,r}}\leq C\|u^2_{n+m}-u^2_{n}\|_{B^{s-1}_{p,r}}\|\partial_x\rho_{n+1}\|_{B^{s-2}_{p,r}}\\
\leq C\|\rho_{n+1}\|_{B^{s-1}_{p,r}}\|u_{n+m}-u_{n}\|_{B^{s-1}_{p,r}}\|u_{n+m}+u_{n}\|_{B^{s-1}_{p,r}}
\leq C \|u_{n+m}-u_{n}\|_{B^{s-1}_{p,r}},
\end{multline}
\begin{align}
\|R^1_{n,m}\|_{B^{s-2}_{p,r}}&\leq \|(\rho_{m+n}-\rho_{n})u_{n+m}\partial_{x}u_{n+m}\|_{B^{s-2}_{p,r}}+\|\rho_{n}(u_{n+m}\partial_{x}u_{n+m}-u_n\partial_xu_n)\|_{B^{s-2}_{p,r}}\\
\nonumber&\leq C(\|\rho_{m+n}-\rho_{n}\|_{B^{s-2}_{p,r}}\|u_{n+m}\partial_{x}u_{n+m}\|_{B^{s-1}_{p,r}}+\|\rho_{n}\|_{B^{s-1}_{p,r}}\|\partial_{x}(u^2_{n+m}-u^2_n)\|_{B^{s-2}_{p,r}})\\
\nonumber&\leq C(\|\rho_{m+n}-\rho_{n}\|_{B^{s-2}_{p,r}}+\|u_{n+m}-u_{n}\|_{B^{s-1}_{p,r}}).
\end{align}
Plugging (3.18) and (3.19) into (3.16) yields that
\begin{multline}
\|\rho_{n+m+1}-\rho_{n+1}\|_{L^{\infty}(0,T;B^{s-2}_{p,r})}\leq C\bigg(\|S_{n+m+1}\rho_0-S_{n+1}\rho_0\|_{B^{s-2}_{p,r}}\\
+\int^T_0\|(\rho_{m+n}-\rho_{n})\|_{B^{s-2}_{p,r}}+\|u_{n+m}-u_{n}\|_{B^{s-1}_{p,r}}dt'\bigg).
\end{multline}
Since $B^{s-1}_{p,r}$ is an algebra, we deduce that
\begin{align}
\|(u^2_{n+m}-u^2_{n})\partial_xu_{n+1}\|_{B^{s-1}_{p,r}}\leq \|u_{n+m}-u_{n}\|_{B^{s-1}_{p,r}}\|u_{n+m}+u_{n}\|_{B^{s-1}_{p,r}}\|\partial_xu_{n+1}\|_{B^{s-1}_{p,r}}\leq C\|u_{n+m}-u_{n}\|_{B^{s-1}_{p,r}}.
\end{align}
 Note that $\|\partial_xG\ast f\|_{B^{s-1}_{p,r}}\leq \|f\|_{B^{s-2}_{p,r}}$ and $\|G\ast f\|_{B^{s-1}_{p,r}}\leq \|f\|_{B^{s-3}_{p,r}}$ , we get
\begin{multline}
\|R^2_{n,m}\|_{B^{s-1}_{p,r}}\leq C(\|u^3_{m+n}-u^3_{n}\|_{B^{s-2}_{p,r}}+\|u_{m+n}(\partial_xu_{m+n})^2-u_{n}(\partial_xu_{n})^2\|_{B^{s-2}_{p,r}}\\
+\|u_{m+n}\rho^2_{m+n}-u_n\rho^2_n\|_{B^{s-2}_{p,r}}+\|(\partial_xu_{m+n})^3-(\partial_xu_{n})^3\|_{B^{s-3}_{p,r}}+\|\partial_xu_{m+n}\rho^2_{m+n}-\partial_xu_n\rho^2_n\|_{B^{s-3}_{p,r}}).
\end{multline}
Since $B^{s-1}_{p,r}$ is an algebra, it follows that
\begin{align}
\|u^3_{m+n}-u^3_{n}\|_{B^{s-2}_{p,r}}\leq \|u^3_{m+n}-u^3_{n}\|_{B^{s-1}_{p,r}}
\leq \|u_{m+n}-u_n\|_{B^{s-1}_{p,r}}\|u^2_{m+n}-u_{m+n}u_n+u^2_n\|_{B^{s-1}_{p,r}}\leq C\|u_{m+n}-u_n\|_{B^{s-1}_{p,r}}.
\end{align}
By virtue of Lemma \ref{Morse}, we infer that
\begin{align}
&\|u_{m+n}(\partial_xu_{m+n})^2-u_{n}(\partial_xu_{n})^2\|_{B^{s-2}_{p,r}}\leq\|(u_{m+n}-u_n)(\partial_xu_{m+n})^2\|_{B^{s-2}_{p,r}}+\|u_n[(\partial_xu_{m+n})^2-(\partial_xu_{n})^2]\|_{B^{s-2}_{p,r}}\\
\nonumber&\leq C(\|u_{m+n}-u_n\|_{B^{s-1}_{p,r}}\|(\partial_xu_{m+n})^2\|_{B^{s-1}_{p,r}}+\|u_n\|_{B^{s-1}_{p,r}}\|(\partial_xu_{m+n})^2-(\partial_xu_{n})^2\|_{B^{s-2}_{p,r}})\\
\nonumber&\leq
C(\|u_{m+n}-u_n\|_{B^{s-1}_{p,r}}+\|\partial_xu_{m+n}-\partial_xu_{n}\|_{B^{s-2}_{p,r}}\|\partial_xu_{m+n}+\partial_xu_{n}\|_{B^{s-1}_{p,r}})\leq C\|u_{m+n}-u_n\|_{B^{s-1}_{p,r}},
\end{align}
\begin{align}
&\|(\partial_xu_{m+n})^3-(\partial_xu_{n})^3\|_{B^{s-3}_{p,r}}\leq\|\partial_x(u_{m+n}-u_n)[(\partial_xu_{m+n})^2-\partial_xu_{m+n}\partial_xu_{n}+(\partial_xu_{n})^2]\|_{B^{s-2}_{p,r}}\\
\nonumber&\leq C(\|\partial_x(u_{m+n}-u_n)\|_{B^{s-2}_{p,r}}\|[(\partial_xu_{m+n})^2-\partial_xu_{m+n}\partial_xu_{n}+(\partial_xu_{n})^2]\|_{B^{s-1}_{p,r}}
\leq C\|u_{m+n}-u_n\|_{B^{s-1}_{p,r}},
\end{align}
\begin{align}
\|u_{m+n}\rho^2_{m+n}-u_n\rho^2_n\|_{B^{s-2}_{p,r}}&\leq \|(u_{m+n}-u_n)\rho^2_{m+n}\|_{B^{s-2}_{p,r}}+\|u_n(\rho^2_{m+n}-\rho^2_n)\|_{B^{s-2}_{p,r}}\\
\nonumber&\leq\|u_{m+n}-u_n\|_{B^{s-2}_{p,r}}\|\rho^2_{m+n}\|_{B^{s-1}_{p,r}}+\|u_n\|_{B^{s-1}_{p,r}}\|\rho^2_{m+n}-\rho^2_n\|_{B^{s-2}_{p,r}}\\
\nonumber&\leq C(\|u_{m+n}-u_n\|_{B^{s-1}_{p,r}}+\|\rho_{m+n}-\rho_n\|_{B^{s-2}_{p,r}}\|\rho_{m+n}+\rho_n\|_{B^{s-1}_{p,r}})\\
\nonumber&\leq C(\|u_{m+n}-u_n\|_{B^{s-1}_{p,r}}+\|\rho_{m+n}-\rho_n\|_{B^{s-2}_{p,r}}),
\end{align}
\begin{align}
\|\partial_xu_{m+n}\rho^2_{m+n}-\partial_xu_n\rho^2_n\|_{B^{s-3}_{p,r}}&\leq \|\partial_x(u_{m+n}-u_n)\rho^2_{m+n}\|_{B^{s-2}_{p,r}}+\|\partial_xu_n(\rho^2_{m+n}-\rho^2_n)\|_{B^{s-2}_{p,r}}\\
\nonumber&\leq\|\partial_x(u_{m+n}-u_n)\|_{B^{s-2}_{p,r}}\|\rho^2_{m+n}\|_{B^{s-1}_{p,r}}+\|\partial_xu_n\|_{B^{s-1}_{p,r}}\|\rho^2_{m+n}-\rho^2_n\|_{B^{s-2}_{p,r}}\\
\nonumber&\leq C(\|u_{m+n}-u_n\|_{B^{s-1}_{p,r}}+\|\rho_{m+n}-\rho_n\|_{B^{s-2}_{p,r}}\|\rho_{m+n}+\rho_n\|_{B^{s-1}_{p,r}})\\
\nonumber&\leq C(\|u_{m+n}-u_n\|_{B^{s-1}_{p,r}}+\|\rho_{m+n}-\rho_n\|_{B^{s-2}_{p,r}}).
\end{align}
Plugging (3.23)-(3.27) into (3.22) yields that
\begin{align}
\|R^2_{n,m}\|_{B^{s-1}_{p,r}}\leq C(\|u_{m+n}-u_n\|_{B^{s-1}_{p,r}}+\|\rho_{m+n}-\rho_n\|_{B^{s-2}_{p,r}}),
\end{align}
which along with (3.17) and (3.21) leads to
\begin{multline}
\|u_{n+m+1}-u_{n+1}\|_{L^{\infty}(0,T;B^{s-1}_{p,r})}\leq C\bigg(\|S_{n+m+1}u_0-S_{n+1}u_0\|_{B^s_{p,r}}\\
+\int^T_0\|u_{m+n}-u_n\|_{B^{s-1}_{p,r}}+\|\rho_{m+n}-\rho_n\|_{B^{s-2}_{p,r}}dt'\bigg).
\end{multline}
The above inequality together with (3.20) yields that
\begin{multline}
\|\rho_{n+m+1}-\rho_{n+1}\|_{L^{\infty}(0,T;B^{s-2}_{p,r})}+\|u_{n+m+1}-u_{n+1}\|_{L^{\infty}(0,T;B^{s-1}_{p,r})}\leq \\ C\bigg(\|S_{n+m+1}\rho_0-S_{n+1}\rho_0\|_{B^{s-1}_{p,r}}+\|S_{n+m+1}u_0-S_{n+1}u_0\|_{B^s_{p,r}}
+\int^T_0\|u_{m+n}-u_n\|_{B^{s-1}_{p,r}}+\|\rho_{m+n}-\rho_n\|_{B^{s-2}_{p,r}}dt'\bigg).
\end{multline}
By defining $A_{n,m}(t)\triangleq \|u_{m+n}-u_n\|_{B^{s-1}_{p,r}}+\|\rho_{m+n}-\rho_n\|_{B^{s-2}_{p,r}}$, we deduce that
\begin{align}
A_{n+1,m}(t)\leq C\bigg(A_{n+1,m}(0)+\int^t_0A_{n,m}(t')dt'\bigg)\leq C\bigg(2^{-n}+\int^t_0A_{n,m}(t')dt'\bigg).
\end{align}
Arguing by induction, we get
\begin{align}
\sup_{t\in[0,T]}A_{n+1,m}(t)\leq \sum^{n}_{k=0}\frac{(CT)^k}{(k+1)!}2^{-n}+\frac{(CT)^{n+1}}{(n+1)!} ,
\end{align}
which leads to $A_{n,m}(t)\rightarrow 0~~as~~n\rightarrow \infty,~\forall m\in\mathbb{N}$. So $\{(\rho_n,u_n)\}$ is a Cauchy sequence in $L^{\infty}(0,T;B^{s-2}_{p,r})\times L^{\infty}(0,T;B^{s-1}_{p,r})$ and converges to some limit function $(\rho,u)\in L^{\infty}(0,T;B^{s-2}_{p,r})\times L^{\infty}(0,T;B^{s-1}_{p,r})$. 
Since $(\rho_n,u_n)$ is bounded in $L^{\infty}(0,T;B^{s-1}_{p,r})\times L^{\infty}(0,T;B^{s}_{p,r})$, the Fatou property for Besov spaces yields that $(\rho,u)$ also belongs to $L^{\infty}(0,T;B^{s-1}_{p,r})\times L^{\infty}(0,T;B^{s}_{p,r})$. An interpolation argument ensures that $(\rho_n,u_n)$ converges to $(\rho,u)$ in $L^{\infty}(0,T;B^{s'-1}_{p,r})\times L^{\infty}(0,T;B^{s'}_{p,r})$ with $s'<s$. Passing to the limit in (3.4), we deduce that $(\rho,u)$ is indeed a solution of (1.2). By virtue of (1.2) and using the fact that $B^s_{p,r}\hookrightarrow C^{0,1}$, we infer that $(\rho,u)\in E^s_{p,r}$. For more details, one can refer to Chapter 3 in \cite{B.C.D}.\\

{\bf Step 4}: Finally, we prove the uniqueness and stability of (1.2). Suppose that $(\rho_1,u_1)$ and $(\rho_2,u_2)$ are two solutions of (1.2). Hence, we obtain that
  \begin{align}
\left\{
\begin{array}{ll}
\partial_t(\rho_{1}-\rho_{2})-u^2_{1}\partial_x(\rho_{1}-\rho_{2})=(u^2_{1}-u^2_{2})\partial_x\rho_{2}+\rho_1u_1\partial_xu_1-\rho_2u_2\partial_xu_2,\\[1ex]
\partial_t(u_{1}-u_{2})-u^2_{1}\partial_x(u_{1}-u_{2})=(u^2_{1}-u^2_{2})\partial_xu_{2}+R(u_1,u_2,\rho_1,\rho_2),\\[1ex]
\end{array}
\right.
\end{align}
where $R(u_1,u_2,\rho_1,\rho_2)=\partial_xG\ast[u^3_{1}-u^3_2+\frac{3}{2}[u_{1}(\partial_xu_{1})^2-u_2(\partial_xu_2)^2]+\frac{1}{2}(u_2\rho^2_2-u_{1}\rho^2_{1})]+\frac{1}{2}G\ast[(\partial_xu_{1})^3-(\partial_xu_2)^3+\partial_xu_2\rho^2_2-\partial_xu_{1}\rho^2_{1}].$
By virtue of Lemma \ref{est1}, we have
\begin{multline}
\|\rho_{1}(t)-\rho_{2}(t)\|_{B^{s-2}_{p,r}}\leq C\bigg(\|\rho_1(0)-\rho_2(0)\|_{B^{s-1}_{p,r}}\\
+\int^t_0\|(u^2_{1}(t')-u^2_{2}(t'))\partial_x\rho_{2}(t')\|_{B^{s-2}_{p,r}}+\|\rho_1(t')u_1(t')\partial_xu_1(t')-\rho_2(t')u_2(t')\partial_xu_2(t')\|_{B^{s-2}_{p,r}}dt'\bigg),
\end{multline}
\begin{multline}
\|u_{1}(t)-u_{2}(t)\|_{B^{s-1}_{p,r}}\leq C\bigg(\|u_1(0)-u_2(0)\|_{B^{s}_{p,r}}\\
+\int^t_0\|(u^2_{1}(t')-u^2_{2}(t'))\partial_xu_{2}(t')\|_{B^{s-1}_{p,r}}+\|R(u_1(t'),u_2(t'),\rho_1(t'),\rho_2(t'))\|_{B^{s-1}_{p,r}}dt'\bigg).
\end{multline}
By a similar calculation as in Step 3, we get
\begin{align}
\|(u^2_{1}-u^2_{2})\partial_x\rho_{2}\|_{B^{s-2}_{p,r}}+\|\rho_1u_1\partial_xu_1-\rho_2u_2\partial_xu_2\|_{B^{s-2}_{p,r}}\leq C(\|\rho_{1}-\rho_{2}\|_{B^{s-2}_{p,r}}+\|u_{1}-u_{2}\|_{B^{s-1}_{p,r}}),
\end{align}
\begin{align}
\|(u^2_{1}-u^2_{2})\partial_xu_{2}\|_{B^{s-1}_{p,r}}\leq C\|u_{1}-u_{2}\|_{B^{s-1}_{p,r}},
\end{align}
\begin{align}
\|R(u_1,u_2,\rho_1,\rho_2)\|_{B^{s-1}_{p,r}}\leq C(\|\rho_{1}-\rho_{2}\|_{B^{s-2}_{p,r}}+\|u_{1}-u_{2}\|_{B^{s-1}_{p,r}}).
\end{align}
Plugging (3.36)-(3.38) into (3.34) and (3.35) yields that
\begin{multline}
\|\rho_{1}(t)-\rho_{2}(t)\|_{B^{s-2}_{p,r}}+\|u_{1}(t)-u_{2}(t)\|_{B^{s-1}_{p,r}}\leq C\bigg(\|\rho_1(0)-\rho_2(0)\|_{B^{s-1}_{p,r}}+\|u_1(0)-u_2(0)\|_{B^{s}_{p,r}}\\
+\int^t_0\|\rho_{1}(t')-\rho_{2}(t')\|_{B^{s-2}_{p,r}}+\|u_{1}(t')-u_{2}(t')\|_{B^{s-1}_{p,r}}dt'\bigg).
\end{multline}
Using Gronwall's inequality, we get
\begin{align}
\sup_{t\in[0,T)}\|\rho_{1}(t)-\rho_{2}(t)\|_{B^{s-2}_{p,r}}+\sup_{t\in[0,T)}\|u_{1}(t)-u_{2}(t)\|_{B^{s-1}_{p,r}}\leq C (\|\rho_1(0)-\rho_2(0)\|_{B^{s-1}_{p,r}}+\|u_1(0)-u_2(0)\|_{B^{s}_{p,r}}).
\end{align}
Taking advantage of the interpolation argument ensures that
\begin{align}
\sup_{t\in[0,T)}\|\rho_{1}(t)-\rho_{2}(t)\|_{B^{s'-1}_{p,r}}+\sup_{t\in[0,T)}\|u_{1}(t)-u_{2}(t)\|_{B^{s'}_{p,r}}\leq C (\|\rho_1(0)-\rho_2(0)\|^{\theta}_{B^{s-1}_{p,r}}+\|u_1(0)-u_2(0)\|^{\theta}_{B^{s}_{p,r}}),
\end{align}
where $\theta=s-s'\in(0,1]$. The above inequality implies the uniqueness. Consequently, we prove the theorem by Steps 1-4.\\
\subsection{Proof of Theorem \ref{th2}}
 Now we turn our attention to prove local well-posedness of (1.2) with initial data in the critical space $B^{\frac{1}{2}}_{2,1}\times B^{\frac{3}{2}}_{2,1}$. \\

 {\bf Step 1}: First, we construct the smooth approximate sequence $(\rho_n, u_n)$ as in the previous subsection. We assume that $(\rho_n, u_n)\in L^{\infty}(0,T;B^{\frac{1}{2}}_{2,1})\times L^{\infty}(0,T;B^{\frac{3}{2}}_{2,1})$. Since $B^{\frac{1}{2}}_{2,1}$ is an algebra, one can check that $\rho_{n}u_n\partial_xu_n\in L^{\infty}(0,T;B^{\frac{1}{2}}_{2,1})$ and $\partial_xG\ast(u^3_n+\frac{3}{2}u_n(\partial_xu_n)^2-u_n\rho_n^2)+\frac{1}{2}G\ast((\partial_xu_n)^3-\partial_xu_n\rho^2_n)\in L^{\infty}(0,T;B^{\frac{3}{2}}_{2,1})$. Taking advantage of the theory of transport equations to (3.4), we obtain that $(\rho_{n+1}, u_{n+1})\in L^{\infty}(0,T;B^{\frac{1}{2}}_{2,1})\times L^{\infty}(0,T;B^{\frac{3}{2}}_{2,1})$. \\

 {\bf Step 2}: By a similar calculation as in Step 2 of the previous subsection . We can find a $T$ satisfies  $4C^3T(\|\rho_0\|_{B^{\frac{1}{2}}_{2,1}}+\|u_0\|_{B^\frac{3}{2}_{2,1}})^2<1$ such that
 \begin{align}
\forall t\in[0,T],~~\|\rho_{n}\|_{B^{\frac{1}{2}}_{2,1}}+\|u_{n}\|_{B^\frac{3}{2}_{2,1}}\leq \frac{C(\|\rho_0\|_{B^{\frac{1}{2}}_{2,1}}+\|u_0\|_{B^s_{p,r}})}{\sqrt{1-4C^3(\|\rho_0\|_{B^{\frac{1}{2}}_{2,1}}+\|u_0\|_{B^\frac{3}{2}_{2,1}})^2 T}},
\end{align}
where $C\geq1$ is a constant independent of $n$ and $T$. \\

{\bf Step 3}: We are going to show that $\{(\rho_n,u_n)\}$ is a Cauchy sequence
 in $L^{\infty}(0,T;B^{-\frac{1}{2}}_{2,\infty})\times L^{\infty}(0,T;B^{\frac{1}{2}}_{2,\infty})$. Applying Lemma \ref{est1} to (3.15), we have
 \begin{multline}
\|\rho_{n+m+1}-\rho_{n+1}\|_{L^{\infty}_t(B^{-\frac{1}{2}}_{2,\infty})}\leq C\bigg(\|S_{n+m+1}\rho_0-S_{n+1}\rho_0\|_{B^{-\frac{1}{2}}_{2,\infty}}\\
+\int^t_0\|(u^2_{n+m}-u^2_{n})\partial_x\rho_{n+1}\|_{B^{-\frac{1}{2}}_{2,\infty}}+\|R^1_{n,m}\|_{B^{-\frac{1}{2}}_{2,\infty}}dt'\bigg),
\end{multline}
\begin{multline}
\|u_{n+m+1}-u_{n+1}\|_{L^{\infty}_t(B^{\frac{1}{2}}_{2,\infty})}\leq C\bigg(\|S_{n+m+1}u_0-S_{n+1}u_0\|_{B^{\frac{1}{2}}_{2,\infty}}\\
+\int^t_0\|(u^2_{n+m}-u^2_{n})\partial_xu_{n+1}\|_{B^{\frac{1}{2}}_{2,\infty}}+\|R^2_{n,m}\|_{B^{\frac{1}{2}}_{2,\infty}}dt'\bigg).
\end{multline}
By virtue of Lemma \ref{Critical}, we deduce that
\begin{multline}
\|(u^2_{n+m}-u^2_{n})\partial_x\rho_{n+1}\|_{B^{-\frac{1}{2}}_{2,\infty}}\leq C\|u^2_{n+m}-u^2_{n}\|_{B^{\frac{1}{2}}_{2,1}}\|\partial_x\rho_{n+1}\|_{B^{-\frac{1}{2}}_{2,\infty}}\\
\leq C\|u_{n+m}-u_{n}\|_{B^{\frac{1}{2}}_{2,1}}\|u_{n+m}+u_{n}\|_{B^{\frac{1}{2}}_{2,1}}\|\rho_{n+1}\|_{B^{\frac{1}{2}}_{2,1}}\leq C\|u_{n+m}-u_{n}\|_{B^{\frac{1}{2}}_{2,1}},
\end{multline}
\begin{align}
\|R^1_{n,m}\|_{B^{-\frac{1}{2}}_{2,\infty}}&\leq \|(\rho_{m+n}-\rho_{n})u_{n+m}\partial_{x}u_{n+m}\|_{B^{-\frac{1}{2}}_{2,\infty}}+\|\rho_{n}(u_{n+m}\partial_{x}u_{n+m}-u_n\partial_xu_n)\|_{B^{-\frac{1}{2}}_{2,\infty}}\\
\nonumber&\leq C(\|\rho_{m+n}-\rho_{n}\|_{B^{-\frac{1}{2}}_{2,\infty}}\|u_{n+m}\partial_{x}u_{n+m}\|_{B^{\frac{1}{2}}_{2,1}}+\|\rho_{n}\|_{B^{\frac{1}{2}}_{2,1}}\|\partial_{x}(u^2_{n+m}-u^2_n)\|_{B^{-\frac{1}{2}}_{2,\infty}})\\
\nonumber&\leq C(\|\rho_{m+n}-\rho_{n}\|_{B^{-\frac{1}{2}}_{2,\infty}}+\|u^2_{n+m}-u^2_{n}\|_{B^{\frac{1}{2}}_{2,1}})\\
\nonumber&\leq C(\|\rho_{m+n}-\rho_{n}\|_{B^{-\frac{1}{2}}_{2,\infty}}+\|u_{n+m}-u_{n}\|_{B^{\frac{1}{2}}_{2,1}}).
\end{align}
Plugging (3.45) and (3.46) into (3.43) yields that
\begin{multline}
\|\rho_{n+m+1}-\rho_{n+1}\|_{L^{\infty}_t(B^{-\frac{1}{2}}_{2,\infty})}\leq C\bigg(\|S_{n+m+1}\rho_0-S_{n+1}\rho_0\|_{B^{-\frac{1}{2}}_{2,\infty}}\\
+\int^t_0\|\rho_{m+n}-\rho_{n}\|_{B^{-\frac{1}{2}}_{2,\infty}}+\|u_{n+m}-u_{n}\|_{B^{\frac{1}{2}}_{2,1}}dt'\bigg).
\end{multline}
Using the fact that $B^{\frac{1}{2}}_{2,1}\hookrightarrow B^{\frac{1}{2}}_{2,\infty}\cap L^\infty$, we obtain
\begin{multline}
\|(u^2_{n+m}-u^2_{n})\partial_xu_{n+1}\|_{B^{\frac{1}{2}}_{2,\infty}}\leq
C\|(u^2_{n+m}-u^2_{n})\partial_xu_{n+1}\|_{B^{\frac{1}{2}}_{2,1}}\leq
C\|(u^2_{n+m}-u^2_{n})\|_{B^{\frac{1}{2}}_{2,1}}\|\partial_xu_{n+1}\|_{B^{\frac{1}{2}}_{2,1}}\\
\leq C\|u_{n+m}-u_{n}\|_{B^{\frac{1}{2}}_{2,1}}\|u_{n+m}+u_{n}\|_{B^{\frac{1}{2}}_{2,1}}\|u_{n+1}\|_{B^{\frac{3}{2}}_{2,1}}\leq C\|u_{n+m}-u_{n}\|_{B^{\frac{1}{2}}_{2,1}}.
\end{multline}
Note that $\|\partial_xG\ast f\|_{B^{\frac{1}{2}}_{2,\infty}}\leq \|f\|_{B^{-\frac{1}{2}}_{2,\infty}}$ and $\|G\ast f\|_{B^{\frac{1}{2}}_{2,\infty}}\leq \|f\|_{B^{-\frac{1}{2}}_{2,\infty}}$. Then, we get
\begin{multline}
\|R^2_{n,m}\|_{B^{\frac{1}{2}}_{2,\infty}}\leq C(\|u^3_{m+n}-u^3_{n}\|_{B^{-\frac{1}{2}}_{2,\infty}}+\|u_{m+n}(\partial_xu_{m+n})^2-u_{n}(\partial_xu_{n})^2\|_{B^{-\frac{1}{2}}_{2,\infty}}\\
+\|u_{m+n}\rho^2_{m+n}-u_n\rho^2_n\|_{B^{-\frac{1}{2}}_{2,\infty}}+\|(\partial_xu_{m+n})^3-(\partial_xu_{n})^3\|_{B^{-\frac{1}{2}}_{2,\infty}}+\|\partial_xu_{m+n}\rho^2_{m+n}-\partial_xu_n\rho^2_n\|_{B^{-\frac{1}{2}}_{2,\infty}}).
\end{multline}
By virtue of Lemma \ref{Critical}, we infer that
\begin{align}
\|u^3_{m+n}-u^3_{n}\|_{B^{-\frac{1}{2}}_{2,\infty}}\leq C\|u_{m+n}-u_n\|_{B^{\frac{1}{2}}_{2,1}}\|u^2_{m+n}-u_{m+n}u_n+u^2_n\|_{B^{-\frac{1}{2}}_{2,\infty}}\leq C\|u_{m+n}-u_n\|_{B^{\frac{1}{2}}_{2,1}},
\end{align}
\begin{align}
&\|u_{m+n}(\partial_xu_{m+n})^2-u_{n}(\partial_xu_{n})^2\|_{B^{-\frac{1}{2}}_{2,\infty}}\leq\|(u_{m+n}-u_n)(\partial_xu_{m+n})^2\|_{B^{-\frac{1}{2}}_{2,\infty}}+\|u_n[(\partial_xu_{m+n})^2-(\partial_xu_{n})^2]\|_{B^{-\frac{1}{2}}_{2,\infty}}\\
\nonumber&\leq C(\|u_{m+n}-u_n\|_{B^{\frac{1}{2}}_{2,1}}\|(\partial_xu_{m+n})^2\|_{B^{-\frac{1}{2}}_{2,\infty}}+\|u_n\|_{B^{\frac{1}{2}}_{2,1}}\|(\partial_xu_{m+n})^2-(\partial_xu_{n})^2\|_{B^{-\frac{1}{2}}_{2,\infty}})\\
\nonumber&\leq
C(\|u_{m+n}-u_n\|_{B^{\frac{1}{2}}_{2,1}}+\|\partial_xu_{m+n}-\partial_xu_{n}\|_{B^{-\frac{1}{2}}_{2,\infty}}\|\partial_xu_{m+n}+\partial_xu_{n}\|_{B^{\frac{1}{2}}_{2,1}})\leq C\|u_{m+n}-u_n\|_{B^{\frac{1}{2}}_{2,1}},
\end{align}
\begin{align}
&\|(\partial_xu_{m+n})^3-(\partial_xu_{n})^3\|_{B^{-\frac{1}{2}}_{2,\infty}}\leq\|\partial_x(u_{m+n}-u_n)[(\partial_xu_{m+n})^2-\partial_xu_{m+n}\partial_xu_{n}+(\partial_xu_{n})^2]\|_{B^{-\frac{1}{2}}_{2,\infty}}\\
\nonumber&\leq C(\|\partial_x(u_{m+n}-u_n)\|_{B^{-\frac{1}{2}}_{2,\infty}}\|[(\partial_xu_{m+n})^2-\partial_xu_{m+n}\partial_xu_{n}+(\partial_xu_{n})^2]\|_{B^{\frac{1}{2}}_{2,1}}
\leq C\|(u_{m+n}-u_n)\|_{B^{\frac{1}{2}}_{2,1}},
\end{align}
\begin{align}
\|u_{m+n}\rho^2_{m+n}-u_n\rho^2_n\|_{B^{-\frac{1}{2}}_{2,\infty}}&\leq \|(u_{m+n}-u_n)\rho^2_{m+n}\|_{B^{-\frac{1}{2}}_{2,\infty}}+\|u_n(\rho^2_{m+n}-\rho^2_n)\|_{B^{-\frac{1}{2}}_{2,\infty}}\\
\nonumber&\leq\|u_{m+n}-u_n\|_{B^{-\frac{1}{2}}_{2,\infty}}\|\rho^2_{m+n}\|_{B^{\frac{1}{2}}_{2,1}}+\|u_n\|_{B^{\frac{1}{2}}_{2,1}}\|\rho^2_{m+n}-\rho^2_n\|_{B^{-\frac{1}{2}}_{2,\infty}}\\
\nonumber&\leq C(\|u_{m+n}-u_n\|_{B^{\frac{1}{2}}_{2,1}}+\|\rho_{m+n}-\rho_n\|_{B^{-\frac{1}{2}}_{2,\infty}}\|\rho_{m+n}+\rho_n\|_{B^{\frac{1}{2}}_{2,1}})\\
\nonumber&\leq C(\|u_{m+n}-u_n\|_{B^{\frac{1}{2}}_{2,1}}+\|\rho_{m+n}-\rho_n\|_{B^{-\frac{1}{2}}_{2,\infty}}),
\end{align}
\begin{align}
\|\partial_xu_{m+n}\rho^2_{m+n}-\partial_xu_n\rho^2_n\|_{B^{-\frac{1}{2}}_{2,\infty}}&\leq \|\partial_x(u_{m+n}-u_n)\rho^2_{m+n}\|_{B^{-\frac{1}{2}}_{2,\infty}}+\|\partial_xu_n(\rho^2_{m+n}-\rho^2_n)\|_{B^{-\frac{1}{2}}_{2,\infty}}\\
\nonumber&\leq\|\partial_x(u_{m+n}-u_n)\|_{B^{-\frac{1}{2}}_{2,\infty}}\|\rho^2_{m+n}\|_{B^{\frac{1}{2}}_{2,1}}+\|\partial_xu_n\|_{B^{\frac{1}{2}}_{2,1}}\|\rho^2_{m+n}-\rho^2_n\|_{B^{-\frac{1}{2}}_{2,\infty}}\\
\nonumber&\leq C(\|u_{m+n}-u_n\|_{B^{\frac{1}{2}}_{2,1}}+\|\rho_{m+n}-\rho_n\|_{B^{-\frac{1}{2}}_{2,\infty}}\|\rho_{m+n}+\rho_n\|_{B^{\frac{1}{2}}_{2,1}})\\
\nonumber&\leq C(\|u_{m+n}-u_n\|_{B^{\frac{1}{2}}_{2,1}}+\|\rho_{m+n}-\rho_n\|_{B^{-\frac{1}{2}}_{2,\infty}}).
\end{align}
Plugging (3.50)-(3.54) into (3.49) yields that
\begin{align}
\|R^2_{n,m}\|_{B^{-\frac{1}{2}}_{2,\infty}}\leq C(\|u_{m+n}-u_n\|_{B^{\frac{1}{2}}_{2,1}}+\|\rho_{m+n}-\rho_n\|_{B^{-\frac{1}{2}}_{2,\infty}}),
\end{align}
which together with  (3.44) and (3.48) leads to
\begin{multline}
\|u_{n+m+1}-u_{n+1}\|_{L^{\infty}_t(B^{\frac{1}{2}}_{2,\infty})}\leq C\bigg(\|S_{n+m+1}u_0-S_{n+1}u_0\|_{B^{\frac{1}{2}}_{2,\infty}}\\
+\int^t_0\|\rho_{m+n}-\rho_{n}\|_{B^{-\frac{1}{2}}_{2,\infty}}+\|u_{n+m}-u_{n}\|_{B^{\frac{1}{2}}_{2,1}}dt'\bigg).
\end{multline}
Hence, we deduce from the above inequality and (3.47) that
\begin{multline}
\|\rho_{m+n+1}-\rho_{n+1}\|_{L^{\infty}_t(B^{-\frac{1}{2}}_{2,\infty})}+\|u_{n+m+1}-u_{n+1}\|_{L^{\infty}_t(B^{\frac{1}{2}}_{2,\infty})}\\
\leq C\bigg(\|S_{n+m+1}\rho_0-S_{n+1}\rho_0\|_{B^{-\frac{1}{2}}_{2,\infty}}+\|S_{n+m+1}u_0-S_{n+1}u_0\|_{B^{\frac{1}{2}}_{2,\infty}}
+\int^t_0\|\rho_{m+n}-\rho_{n}\|_{B^{-\frac{1}{2}}_{2,\infty}}+\|u_{n+m}-u_{n}\|_{B^{\frac{1}{2}}_{2,1}}dt'\bigg).
\end{multline}
Applying Lemma \ref{log.est} to the above inequality, we have
 \begin{multline}
\|\rho_{m+n+1}-\rho_{n+1}\|_{L^{\infty}_t(B^{-\frac{1}{2}}_{2,\infty})}+\|u_{n+m+1}-u_{n+1}\|_{L^{\infty}_t(B^{\frac{1}{2}}_{2,\infty})}\\
\leq C\bigg(a_{n,m}
+\int^t_0\|\rho_{m+n}-\rho_{n}\|_{B^{-\frac{1}{2}}_{2,\infty}}+\|u_{n+m}-u_{n}\|_{B^{\frac{1}{2}}_{2,\infty}}\ln(e+\frac{\|u_{n+m}-u_{n}\|_{B^{\frac{3}{2}}_{2,1}}}{\|u_{n+m}-u_{n}\|_{B^{\frac{1}{2}}_{2,1}}})dt'\bigg),
\end{multline}
where $a_{n,m}=\|S_{n+m+1}\rho_0-S_{n+1}\rho_0\|_{B^{-\frac{1}{2}}_{2,\infty}}+\|S_{n+m+1}u_0-S_{n+1}u_0\|_{B^{\frac{1}{2}}_{2,\infty}}$.
Since the function $x\ln(e+\frac{c}{x})$ with $c>0$ is nondecreasing, it follows that
 \begin{multline}
\|\rho_{m+n+1}-\rho_{n+1}\|_{L^{\infty}_t(B^{-\frac{1}{2}}_{2,\infty})}+\|u_{n+m+1}-u_{n+1}\|_{L^{\infty}_t(B^{\frac{1}{2}}_{2,\infty})}\\
\leq C\bigg(a_{n,m}
+\int^t_0(\|\rho_{m+n}-\rho_{n}\|_{B^{-\frac{1}{2}}_{2,\infty}}+\|u_{n+m}-u_{n}\|_{B^{\frac{1}{2}}_{2,\infty}})\ln(e+\frac{C}{(\|\rho_{m+n}-\rho_{n}\|_{B^{-\frac{1}{2}}_{2,\infty}}+\|u_{n+m}-u_{n}\|_{B^{\frac{1}{2}}_{2,\infty}})})dt'\bigg),
\end{multline}
By defining $A_{n,m}(t)\triangleq \|\rho_{m+n}-\rho_{n}\|_{B^{-\frac{1}{2}}_{2,\infty}}+\|u_{n+m}-u_{n}\|_{B^{\frac{1}{2}}_{2,\infty}}$, we obtain
\begin{align}
A_{n+1,m}(t)\leq C\bigg(a_{n,m}+\int^t_0A_{n,m}(t')\ln(e+\frac{C}{A_{n,m}(t')})dt'\bigg)\leq C\bigg(a_{n,m}+\int^t_0A_{n,m}(t')(1-\ln\frac{A_{n,m}(t')}{C})dt'\bigg).
\end{align}
Let $A_n(t)=\sup_{m}A_{n,m}(t)$ and $a_n=\sup_m a_{n,m}$. As the function $x(1-\ln{\frac{x}{C}})$ is nondecreasing in $x\in [0,C)$, we get
\begin{align}
A_{n+1}(t)\leq C\bigg(a_n+\int^t_0A_n(t')(1-\ln\frac{A_n(t')}{C})dt'\bigg),
\end{align}
which together with Lebesgue's dominated convergence theorem leads to
\begin{align}
\widetilde{A}(t)\triangleq \limsup_{n\rightarrow\infty}A_{n+1}(t)\leq C\int^t_0\widetilde{A}(t')(1-\ln\frac{\widetilde{A}(t')}{C})dt.
\end{align}
Taking advantage of Lemma \ref{Osgood}, we deduce that $\widetilde{A}(t)=0$. In other words  $\{(\rho_n,u_n)\}$ is a Cauchy sequence
 in $L^{\infty}(0,T;B^{-\frac{1}{2}}_{2,\infty})\times L^{\infty}(0,T;B^{\frac{1}{2}}_{2,\infty})$ and converges to some limit function $(\rho,u)\in L^{\infty}(0,T;B^{-\frac{1}{2}}_{2,\infty})\times L^{\infty}(0,T;B^{\frac{1}{2}}_{2,\infty})$ . By the similar argument as in Step 3 of the previous subsection, we verify that $(\rho,u)\in C([0,T);B^{\frac{1}{2}}_{2,1})\times C([0,T);B^{\frac{3}{2}}_{2,1})\cap C^1([0,T);B^{-\frac{1}{2}}_{2,1})\times C^1([0,T);B^{\frac{1}{2}}_{2,1})$ is indeed a solution of (1.2).\\

 {\bf Step 4}: Finally, we prove the uniqueness of (1.2).  Applying Lemma \ref{est1} to (3.33), we have
 \begin{multline}
\|\rho_{1}(t)-\rho_{2}(t)\|_{B^{-\frac{1}{2}}_{2,\infty}}\leq C\bigg(\|\rho_1(0)-\rho_2(0)\|_{B^{-\frac{1}{2}}_{2,\infty}}\\
+\int^t_0\|(u^2_{1}(t')-u^2_{2}(t'))\partial_x\rho_{2}(t')\|_{B^{-\frac{1}{2}}_{2,\infty}}+\|\rho_1(t')u_1(t')\partial_xu_1(t')-\rho_2(t')u_2(t')\partial_xu_2(t')\|_{B^{-\frac{1}{2}}_{2,\infty}}dt'\bigg),
\end{multline}
\begin{multline}
\|u_{1}(t)-u_{2}(t)\|_{B^{\frac{1}{2}}_{2,\infty}}\leq C\bigg(\|u_1(0)-u_2(0)\|_{B^{\frac{1}{2}}_{2,\infty}}\\
+\int^t_0\|(u^2_{1}(t')-u^2_{2}(t'))\partial_xu_{2}(t')\|_{B^{\frac{1}{2}}_{2,\infty}}+\|R(u_1(t'),u_2(t'),\rho_1(t'),\rho_2(t'))\|_{B^{\frac{1}{2}}_{2,\infty}}dt'\bigg).
\end{multline}
By a similar calculation as in Step 3, we get
\begin{align}
\|(u^2_{1}-u^2_{2})\partial_x\rho_{2}\|_{B^{-\frac{1}{2}}_{2,\infty}}+\|\rho_1u_1\partial_xu_1-\rho_2u_2\partial_xu_2\|_{B^{-\frac{1}{2}}_{2,\infty}}\leq C(\|\rho_{1}-\rho_{2}\|_{B^{-\frac{1}{2}}_{2,\infty}}+\|u_{1}-u_{2}\|_{B^{\frac{1}{2}}_{2,1}}),
\end{align}
\begin{align}
\|(u^2_{1}-u^2_{2})\partial_xu_{2}\|_{B^{\frac{1}{2}}_{2,\infty}}\leq C\|u_{1}-u_{2}\|_{B^{\frac{1}{2}}_{2,1}},
\end{align}
\begin{align}
\|R(u_1,u_2,\rho_1,\rho_2)\|_{B^{\frac{1}{2}}_{2,\infty}}\leq C(\|\rho_{1}-\rho_{2}\|_{B^{-\frac{1}{2}}_{2,\infty}}+\|u_{1}-u_{2}\|_{B^{\frac{1}{2}}_{2,1}}).
\end{align}
Plugging (3.65)-(3.67) into (3.63) and (3.64) yields that
\begin{multline}
\|\rho_{1}(t)-\rho_{2}(t)\|_{B^{-\frac{1}{2}}_{2,\infty}}+\|u_{1}(t)-u_{2}(t)\|_{B^{\frac{1}{2}}_{2,\infty}}\leq C\bigg(\|\rho_1(0)-\rho_2(0)\|_{B^{-\frac{1}{2}}_{2,\infty}}+\|u_1(0)-u_2(0)\|_{B^{\frac{1}{2}}_{2,\infty}}\\
+\int^t_0\|\rho_{1}(t')-\rho_{2}(t')\|_{B^{-\frac{1}{2}}_{2,\infty}}+\|u_{1}(t')-u_{2}(t')\|_{B^{\frac{1}{2}}_{2,1}}dt'\bigg).
\end{multline}
Applying Lemma \ref{log.est} to the above inequality, we have
\begin{multline}
\|\rho_{1}(t)-\rho_{2}(t)\|_{B^{-\frac{1}{2}}_{2,\infty}}+\|u_{1}(t)-u_{2}(t)\|_{B^{\frac{1}{2}}_{2,\infty}}\leq C\bigg(\|\rho_1(0)-\rho_2(0)\|_{B^{-\frac{1}{2}}_{2,\infty}}+\|u_1(0)-u_2(0)\|_{B^{\frac{1}{2}}_{2,\infty}}\\
+\int^t_0\|\rho_{1}(t')-\rho_{2}(t')\|_{B^{-\frac{1}{2}}_{2,\infty}}+\|u_{1}(t')-u_{2}(t')\|_{B^{\frac{1}{2}}_{2,\infty}}\ln(e+\frac{C}{\|u_{1}(t')-u_{2}(t')\|_{B^{\frac{1}{2}}_{2,\infty}}})dt'\bigg).
\end{multline}
Since the function $x\ln(e+\frac{c}{x})$ is nondecreasing, it follows that
 \begin{align}
A(t)\leq C\bigg(A(0)
+\int^t_0(A(t')\ln(e+\frac{C}{A(t')})dt'\bigg)\leq C\bigg(A(0)
+\int^t_0A(t')(1-\ln\frac{A(t')}{C})dt'\bigg),
\end{align}
where $A(t)=\|\rho_{1}(t)-\rho_{2}(t)\|_{B^{-\frac{1}{2}}_{2,\infty}}+\|u_{1}(t)-u_{2}(t)\|_{B^{\frac{1}{2}}_{2,\infty}}$. By virtue of Lemma \ref{Osgood}, we verify that
\begin{align}
\nonumber\|\rho_{1}(t)-\rho_{2}(t)\|_{B^{-\frac{1}{2}}_{2,\infty}}+\|u_{1}(t)-u_{2}(t)\|_{B^{\frac{1}{2}}_{2,\infty}}&\leq C(\|\rho_{1}(0)-\rho_{2}(0)\|_{B^{-\frac{1}{2}}_{2,\infty}}+\|u_{1}(0)-u_{2}(0)\|_{B^{\frac{1}{2}}_{2,\infty}})^{\exp\{-Ct\}}\\
&\leq C(\|\rho_{1}(0)-\rho_{2}(0)\|_{B^{\frac{1}{2}}_{2,1}}+\|u_{1}(0)-u_{2}(0)\|_{B^{\frac{3}{2}}_{2,1}})^{\exp\{-Ct\}},
\end{align}
Taking advantage of the interpolation argument ensures that
\begin{multline}
\sup_{t\in[0,T)}\|\rho_{1}(t)-\rho_{2}(t)\|_{B^{s'-1}_{2,1}}+\sup_{t\in[0,T)}\|u_{1}(t)-u_{2}(t)\|_{B^{s'}_{2,1}}\\
\leq C (\|\rho_1(0)-\rho_2(0)\|_{B^{\frac{1}{2}}_{2,1}}+\|u_1(0)-u_2(0)\|_{B^{\frac{3}{2}}_{2,1}})^{\theta\exp\{-Ct\}},
\end{multline}
where $\theta=\frac{3}{2}-s'\in(0,1]$. The above inequality implies the uniqueness. Consequently, we prove the theorem by Steps 1-4.\\
\section{Blow-up criteria}
In this section, we present two blow-up criteria for (1.1). Our first result can be stated as follows.
\begin{theo}\label{th3}
Let $(\rho_0,u_0)\in B^{\frac{1}{2}}_{2,1}(\mathbb{R})\times B^{\frac{3}{2}}_{2,1}(\mathbb{R})$ and let $T$ be the maximal existence time of the solution $(\rho,u)$ to (1.2). Then the solution blows up in finite time if and only if  $$\int^T_0\|u(t)\|^2_{L^\infty}+\|u_x(t)\|^2_{L^\infty}+\|\rho\|^2_{L^\infty}dt=\infty.$$
\begin{proof}
Applying $\Delta_j$ to (1.2) yields that
  \begin{align}
\left\{
\begin{array}{ll}
\Delta_j\rho_{t}-\Delta_j\rho_{x}u^2=R^1_j+\Delta_j(\rho uu_{x}),\\[1ex]
\Delta_ju_{t}-\Delta_ju_xu^2=R^2_j+\Delta_j\partial_xG\ast(u^3+\frac{3}{2}uu^2_x-\frac{1}{2}u\rho^2)+\frac{1}{2}\Delta_jG\ast(u^3_x-u_x\rho^2) ,\\[1ex]
\Delta_j\rho|_{t=0}=\Delta_j\rho_{0}, \Delta_ju|_{t=0}=\Delta_ju_{0},\\[1ex]
\end{array}
\right.
\end{align}
where $R^1_j=\Delta_j(\rho_{x}u^2)-\Delta_j\rho_{x}u^2$ and $R^2_j=\Delta_j(u_{x}u^2)-\Delta_ju_{x}u^2$. Multiplying both sides of the first equation of (4.1) by $\Delta_j\rho$ and integrating over $\mathbb{R}$ with respect to $x$, we obtain
\begin{align}
\frac{d}{dt}\frac{1}{2}\int_{\mathbb{R}}(\Delta_j\rho)^2-\int_{\mathbb{R}}\Delta_j\rho\Delta_j\rho_{x}u^2=\int_{\mathbb{R}}\Delta_j\rho R^1_j+\int_{\mathbb{R}}\Delta_j\rho\Delta_j(\rho uu_{x}).
\end{align}
Using H\"{o}lder's inequality and integration by parts, we deduce that
\begin{align}
\|\Delta_j\rho\|_{L^2}\leq \|\Delta_j\rho_0\|_{L^2}+C\displaystyle\int^t_0\|\Delta_j\rho\|_{L^2}\|u\|_{L^\infty}\|u_x\|_{L^\infty}+\|R^1_j\|_{L^2}+\|\Delta_j(\rho uu_{x})\|_{L^2}dt'.
\end{align}
By the same token, we get
\begin{multline}
\|\Delta_ju\|_{L^2}\leq \|\Delta_ju_0\|_{L^2}+C\displaystyle\int^t_0\|\Delta_ju\|_{L^2}\|u\|_{L^\infty}\|u_x\|_{L^\infty}+\|R^2_j\|_{L^2}\\
+\|\Delta_j\partial_xG\ast(u^3+\frac{3}{2}uu^2_x-\frac{1}{2}u\rho^2)+\frac{1}{2}\Delta_jG\ast(u^3_x-u_x\rho^2)\|_{L^2}dt'.
\end{multline}
Multiplying both sides of (4.3) by $2^{\frac{1}{2}j}$ and taking $l^1$-norm, we deduce that
\begin{align}
\|\rho\|_{B^{\frac{1}{2}}_{2,1}}\leq \|\rho_0\|_{B^{\frac{1}{2}}_{2,1}}+C\displaystyle\int^t_0\|\rho\|_{B^{\frac{1}{2}}_{2,1}}\|u\|_{L^\infty}\|u_x\|_{L^\infty}+\|2^{\frac{1}{2}j}\|R^1_j\|_{L^2}\|_{l^1}+\|\rho uu_{x}\|_{B^{\frac{1}{2}}_{2,1}}dt'.
\end{align}
By virtue of Lemma \ref{Rj.est}, we obtain
\begin{align}
\|2^{\frac{1}{2}j}\|R^1_j\|_{L^2}\|_{l^1} \leq C\|\rho\|_{B^{\frac{1}{2}}_{2,1}}\|(u^2)_x\|_{L^\infty} \leq C\|\rho\|_{B^{\frac{1}{2}}_{2,1}}\|u\|_{L^\infty}\|u_x\|_{L^\infty}.
\end{align}
Plugging (4.6) into (4.5) yields that
\begin{align}
\|\rho\|_{B^{\frac{1}{2}}_{2,1}}&\leq \|\rho_0\|_{B^{\frac{1}{2}}_{2,1}}+C\displaystyle\int^t_0\|\rho\|_{B^{\frac{1}{2}}_{2,1}}\|u\|_{L^\infty}\|u_x\|_{L^\infty}+\|\rho uu_{x}\|_{B^{\frac{1}{2}}_{2,1}}dt'\\
\nonumber&\leq \|\rho_0\|_{B^{\frac{1}{2}}_{2,1}}+C\displaystyle\int^t_0\|\rho\|_{B^{\frac{1}{2}}_{2,1}}\|u\|_{L^\infty}\|u_x\|_{L^\infty}+\|\rho\|_{L^\infty}\|(u^2)_x\|_{B^{\frac{1}{2}}_{2,1}}dt'\\
\nonumber&\leq
\|\rho_0\|_{B^{\frac{1}{2}}_{2,1}}+C\displaystyle\int^t_0\|\rho\|_{B^{\frac{1}{2}}_{2,1}}\|u\|_{L^\infty}\|u_x\|_{L^\infty}+\|\rho\|_{L^\infty}\|u\|_{L^\infty}\|u\|_{B^{\frac{3}{2}}_{2,1}}dt'.
\end{align}
Multiplying both sides of (4.4) by $2^{\frac{3}{2}j}$ and taking $l^1$-norm, we infer that
\begin{multline}
\|u\|_{B^{\frac{3}{2}}_{2,1}} \leq \|u_0\|_{B^{\frac{3}{2}}_{2,1}}+C\displaystyle\int^t_0\|u\|_{B^{\frac{3}{2}}_{2,1}}\|u\|_{L^\infty}\|u_x\|_{L^\infty}+\|2^{\frac{3}{2}j}\|R^2_j\|_{L^2}\|_{l^1}\\
+\|u^3\|_{B^{\frac{1}{2}}_{2,1}}+\|uu^2_x\|_{B^{\frac{1}{2}}_{2,1}}+\|u\rho^2\|_{B^{\frac{1}{2}}_{2,1}}+\|u^3_x\|_{B^{\frac{1}{2}}_{2,1}}+\|u_x\rho^2\|_{B^{\frac{1}{2}}_{2,1}}dt'.
\end{multline}
By virtue of Lemma \ref{Rj.est}, we get
\begin{align}
\|2^{\frac{3}{2}j}\|R^2_j\|_{L^2}\|_{l^1} \leq C\|u\|_{B^{\frac{3}{2}}_{2,1}}\|(u^2)_x\|_{L^\infty} \leq C\|u\|_{B^{\frac{3}{2}}_{2,1}}\|u\|_{L^\infty}\|u_x\|_{L^\infty}.
\end{align}
Plugging (4.9) into (4.8) yields that
\begin{align}
\|u\|_{B^{\frac{3}{2}}_{2,1}} &\leq \|u_0\|_{B^{\frac{3}{2}}_{2,1}}+C\displaystyle\int^t_0\|u\|_{B^{\frac{3}{2}}_{2,1}}\|u\|_{L^\infty}\|u_x\|_{L^\infty}
\\
\nonumber&+\|u^3\|_{B^{\frac{1}{2}}_{2,1}}+\|uu^2_x\|_{B^{\frac{1}{2}}_{2,1}}+\|u\rho^2\|_{B^{\frac{1}{2}}_{2,1}}+\|u^3_x\|_{B^{\frac{1}{2}}_{2,1}}+\|u_x\rho^2\|_{B^{\frac{1}{2}}_{2,1}}dt'\\
\nonumber&\leq \|u_0\|_{B^{\frac{3}{2}}_{2,1}}+C\displaystyle\int^t_0 \|u\|_{B^{\frac{3}{2}}_{2,1}}(\|u\|^2_{L^\infty}+\|u_x\|^2_{L^\infty}+\|\rho\|^2_{L^\infty})\\
\nonumber&+\|\rho\|_{B^{\frac{1}{2}}_{2,1}}\|\rho\|_{L^\infty}(\|u\|_{L^\infty}+\|u_x\|_{L^\infty})dt'.
\end{align}
Combining (4.7) and (4.8), we have
\begin{multline}
\|\rho\|_{B^{\frac{1}{2}}_{2,1}}+\|u\|_{B^{\frac{3}{2}}_{2,1}}\leq \|\rho_0\|_{B^{\frac{1}{2}}_{2,1}}+\|u_0\|_{B^{\frac{3}{2}}_{2,1}}+
C\displaystyle\int^t_0(\|u\|_{B^{\frac{3}{2}}_{2,1}}+\|\rho\|_{B^{\frac{1}{2}}_{2,1}})(\|u\|^2_{L^\infty}+\|u_x\|^2_{L^\infty}+\|\rho\|^2_{L^\infty})dt'.
\end{multline}
Taking advantage of Gronwall's inequality, we get
\begin{align}
\|\rho\|_{B^{\frac{1}{2}}_{2,1}}+\|u\|_{B^{\frac{3}{2}}_{2,1}}\leq (\|\rho_0\|_{B^{\frac{1}{2}}_{2,1}}+\|u_0\|_{B^{\frac{3}{2}}_{2,1}})\exp\{\int^t_0C(\|u\|^2_{L^\infty}+\|u_x\|^2_{L^\infty}+\|\rho\|^2_{L^\infty})dt'\}.
\end{align}
Therefore, if $T<\infty$ satisfies that  $\int^T_0\|u\|^2_{L^\infty}+\|u_x\|^2_{L^\infty}+\|\rho\|^2_{L^\infty}dt'<\infty$, then we deduce from (4.12) that
\begin{align}
 \limsup_{t\rightarrow T}(\|\rho(t)\|_{B^{\frac{1}{2}}_{2,1}}+\|u(t)\|_{B^{\frac{3}{2}}_{2,1}})<\infty,
 \end{align}
 which contradicts the assumption that $T$ is the maximal existence time.
\end{proof}
\end{theo}
\begin{rema}
Following the similar proof of Theorem \ref{th3}, one can obtain the same blow-up criterion for (1.2) with initial data $(\rho_0, u_0)$ satisfies the condition of Theorem \ref{th1}.
\end{rema}
In order to obtain the second criterion, we need to draw a support from the following ordinary differential equation:
   \begin{align}
\left\{
\begin{array}{ll}
\frac{d\Phi(t,x)}{dt}=-u^2(t,\Phi(t,x)),\\[1ex]
\Phi(0,x)=x,
\end{array}
\right.
\end{align}
 for the flow generated by $-u^2(t,x)$. Since $u\in C([0,T);B^{\frac{3}{2}}_{2,1})\hookrightarrow C([0,T);Lip)$, it follows that $u^2\in C([0,T);Lip)$.
 Applying classical results in the theory of ordinary differential equations, one can obtain the following lemma.
 \begin{lemm}\cite{LY.YZY}\label{Flow}
 Let  $(\rho_0,u_0)\in B^{\frac{1}{2}}_{2,1}(\mathbb{R})\times B^{\frac{3}{2}}_{2,1}(\mathbb{R})$ and let $T > 0$ be the maximal existence time of the corresponding solution $(\rho,u)$ to (1.2). Then (4.14) has a unique solution
$\Phi \in C^1([0, T )\times\mathbb{R};\mathbb{R})$. Moreover, the map $\Phi(t, \cdot)$ is an increasing diffeomorphism
of $\mathbb{R}$ with
\begin{align}
\Phi_x(t,x)=\exp\bigg(\int^t_0-2uu_x(s,\Phi(s,x))ds\bigg)>0,~~ \forall (t,x)\in[0,T)\times\mathbb{R}.
\end{align}
 \end{lemm}
 Now we introduce two conservation laws for (1.2) which are crucial to obtain our result.
 \begin{lemm}\label{conservation}
 Let $(\rho_0,u_0)\in B^{\frac{1}{2}}_{2,1}(\mathbb{R})\times B^{\frac{3}{2}}_{2,1}(\mathbb{R})$ and let $T$ be the maximal existence time of the solution $(\rho,u)$ to (1.2). Then for all $t\in [0,T)$, we have
 \begin{align}
 \int_{\mathbb{R}}u^2(t)+u^2_x(t)dx=\int_{\mathbb{R}}u^2_0+u^2_{0,x}dx,~~~~~\int_{\mathbb{R}}\rho^2(t)dx=\int_{\mathbb{R}}\rho_0^2dx.
 \end{align}
 Moreover, $\|u(t)\|_{L^\infty}\leq C\|u(t)\|_{H^1}=C\|u_0\|_{H^1}$.
 \begin{proof}
  Arguing by density, it suffices to consider the case where $(\rho, u)\in C^{\infty}_0(\mathbb{R})$. Taking advantage of (1.1) and integration by part, we deduce that
  \begin{align}
  \frac{d}{dt}\int_{\mathbb{R}}u^2+u^2_xdx&=2\int_{\mathbb{R}}um_tdx=2\int_{\mathbb{R}}3u_xu^2m+u^3m_x-u\rho(u\rho)_xdx\\
  \nonumber&=\int_{\mathbb{R}}2(u^3m)_x-[(\rho u)^2]_xdx=0
  \end{align}
  and
  \begin{align}
  \frac{d}{dt}\int_{\mathbb{R}}\rho^2dx=2\int_{\mathbb{R}}\rho\rho_tdx=\int_{\mathbb{R}}[(\rho u)^2]_xdx=0.
  \end{align}
 \end{proof}
 \end{lemm}

Thanks to the above lemmas we have the following result.
 \begin{theo}\label{th4}
 Let $(\rho_0,u_0)\in B^{\frac{1}{2}}_{2,1}(\mathbb{R})\times B^{\frac{3}{2}}_{2,1}(\mathbb{R})$ and let $T$ be the maximal existence time of the solution $(\rho,u)$ to (1.2). Then the solution blows up in finite time if and only if
 $$\limsup_{t\rightarrow T}\sup_{x\in\mathbb{R}}uu_x(t,x)=+\infty.$$
 \begin{proof}
 By virtue of (1.2), we infer that
 \begin{align}
 \frac{d}{dt}\rho(t,\Phi(t,x))=\rho uu_x(t,\Phi(t,x)),
 \end{align}
 which leads to
 \begin{align}
 \rho(t,\Phi(t,x))=\rho_0 e^{\int^t_0uu_x(s,\Phi(s,x)))ds}.
 \end{align}
  Assume that the solution $(\rho,u)$ of (1.2) blows up in finite time $(T<\infty)$ and there exists a constant $M>0$ such that
  \begin{align}
  \sup_{x\in\mathbb{R}}uu_x(t,x)\leq M.
  \end{align}
  Thanks to Lemma \ref{Flow}, we obtain
  \begin{align}
  \|\rho(t,\cdot)\|_{L^\infty}=\|\rho(t,\Phi(t,\cdot))\|_{L^\infty}\leq \|\rho_0\|_{L^\infty}e^{MT}<\infty.
  \end{align}
  By differentiating the second equation of (1.2) with respect to $x$, we deduce that
  \begin{align}
  u_{xt}=u^2u_{xx}+\frac{1}{2}uu^2_x-u^3+u\rho^2+G\ast(u^3+\frac{3}{2}uu^2_x-\frac{1}{2}u\rho^2)+\frac{1}{2}\partial_xG\ast(u^3_x-u_x\rho^2),
  \end{align}
 which leads to
 \begin{multline}
 \frac{d}{dt}u_x(t,\Phi(t,x))=\frac{1}{2}uu^2_x(t,\Phi(t,x))-u^3(t,\Phi(t,x))+u\rho^2(t,\Phi(t,x))\\
 +G\ast(u^3+\frac{3}{2}uu^2_x-\frac{1}{2}u\rho^2)(t,\Phi(t,x))+\frac{1}{2}\partial_xG\ast(u^3_x-u_x\rho^2)(t,\Phi(t,x)).
 \end{multline}
 Integrating over $[0,t]$ with respect to $t$, we get
\begin{multline}
u_x(t,\Phi(t,x))=u_{0,x}+\int^t_0\frac{1}{2}uu^2_x(s,\Phi(s,x))-u^3(s,\Phi(s,x))+u\rho^2(s,\Phi(s,x))\\
+G\ast(u^3+\frac{3}{2}uu^2_x-\frac{1}{2}u\rho^2)(s,\Phi(s,x))+\frac{1}{2}\partial_xG\ast(u^3_x-u_x\rho^2)(s,\Phi(s,x))ds.
\end{multline}
By virtue of Lemma \ref{conservation}, we see that
\begin{align}
G\ast(u^3+\frac{3}{2}uu^2_x-\frac{1}{2}u\rho^2)(t,\Phi(t,x))&=\frac{1}{2}\int_{\mathbb{R}}e^{-|\Phi(t,x)-y|}(u^3+\frac{3}{2}uu^2_y-\frac{1}{2}u\rho^2)(t,y)dy\\
\nonumber&\leq \|u\|_{L^\infty}(\|u\|_{H^1}+\|\rho\|_{L^2})\leq C\|u_0\|_{H^1}(\|u_0\|_{H^1}+\|\rho_0\|_{L^2}),
\end{align}
\begin{align}
\partial_xG\ast(u^3_x-u_x\rho^2)(s,\Phi(s,x))\leq \|u_x\|_{L^\infty}(\|u\|_{H^1}+\|\rho\|_{L^2})\leq \|u_x\|_{L^\infty}(\|u_0\|_{H^1}+\|\rho_0\|_{L^2}).
\end{align}
Plugging (4.21), (4.22), (4.26) and (4.27) into (4.25) yields that
\begin{align}
u_x(t,\Phi(t,x))\leq u_{0,x}+\int^t_0(\frac{M}{2}+\|u_0\|_{H^1}+\|\rho_0\|_{L^2})\|u_x\|_{L^\infty}+C\|u_0\|_{H^1}[\|u_0\|_{H^1}+\|\rho_0\|_{L^2}(1+e^{MT})]ds.
\end{align}
Taking advantage of Gronwall's inequality, we deduce that
\begin{align}
\|u_x\|_{L^\infty}< \bigg\{\|u_{0,x}\|_{L^\infty}+CT\|u_0\|_{H^1}[\|u_0\|_{H^1}+\|\rho_0\|_{L^2}(1+e^{MT})]\bigg\}e^{(\frac{M}{2}+\|u_0\|_{H^1}+\|\rho_0\|_{L^2})T}<\infty,
\end{align}
which contradicts to Theorem \ref{th3}.\\
 On the other hand, by Theorem \ref{th3} and Sobolev's embedding theorem, we see that if
 $$\limsup_{t\rightarrow T}\sup_{x\in\mathbb{R}}uu_x(t,x)=+\infty$$
 then the solution $(\rho,u)$ will blow up in finite time. This completes the proof of the theorem.
 \end{proof}
 \end{theo}
 \begin{rema}
 If $(\rho_0,u_0)\in L^2\times H^1$ and satisfies the condition of Theorem \ref{th1}, one can also obtain the same blow-up criterion for (1.2).
 \end{rema}
 \begin{rema}
 If $\rho=0$, then Theorem \ref{th4} covers the recent results for the Novikov equation in \cite{Wu.Yin2,Wei.Yan2}.
 \end{rema}

 {\bf Acknowledgements}.  This work was
partially supported by NNSFC (No.11271382), RFDP (No.
20120171110014), the Macau Science and Technology Development Fund (No. 098/2013/A3) and the key project of Sun Yat-sen University. The authors thank the referees for their valuable comments and suggestions.

\end{document}